\documentclass[11pt]{amsart}
\usepackage{fullpage}
\usepackage{color}
\usepackage{pstricks,pst-node,pst-plot} 
\usepackage{graphicx,psfrag} 
\usepackage{color} 
\usepackage{tikz}
\usepackage{pgffor}

\usepackage{perpage}	 

\MakePerPage{footnote}	 

\newtheorem{theorem}{Theorem}[section]
\newtheorem{lemma}[theorem]{Lemma}
\newtheorem{corollary}[theorem]{Corollary}
\newtheorem{proposition}[theorem]{Proposition}

\theoremstyle{definition}
\newtheorem{definition}[theorem]{Definition}
\newtheorem{example}[theorem]{Example}

\newtheorem{question}[theorem]{Question}

\theoremstyle{remark}
\newtheorem{remark}[theorem]{Remark}

\numberwithin{equation}{section}

\definecolor{gray}{rgb}{.5,.5,.5}

\definecolor{black}{rgb}{0,0,0}

\definecolor{blue}{rgb}{0,0,1}

\definecolor{red}{rgb}{1,0,0}

\definecolor{green}{rgb}{0,1,0}

\definecolor{yellow}{rgb}{1,1,.4}

\newrgbcolor{purple}{.5 0 .5}
\newrgbcolor{black}{0 0 0}
\newrgbcolor{white}{1 1 1}
\newrgbcolor{gold}{.5 .5 .2}
\newrgbcolor{darkgreen}{0 .5 0}
\newrgbcolor{gray}{.5 .5 .5}
\newrgbcolor{lightgray}{.75 .75 .75}

\begin{document}

\markboth{Moshe Cohen}
{A determinant formula for the Jones polynomial of pretzel knots}

\title{A determinant formula for\\ the Jones polynomial of pretzel knots}



\author{Moshe Cohen}

\address{Department of Mathematics and Computer Science, Bar-Ilan University, Ramat Gan 52900, Israel}

\email{cohenm10@macs.biu.ac.il}

\urladdr{http://u.math.biu.ac.il/~cohenm10}

\begin{abstract}
This paper presents an algorithm to construct a weighted adjacency matrix of a plane bipartite graph obtained from a pretzel knot diagram.  The determinant of this matrix after evaluation is shown to be the Jones polynomial of the pretzel knot by way of perfect matchings (or dimers) of this graph.  The weights are Tutte's activity letters that arise because the Jones polynomial is a specialization of the signed version of the Tutte polynomial.  The relationship is formalized between the familiar spanning tree setting for the Tait graph and the perfect matchings of the plane bipartite graph above.  Evaluations of these activity words are related to the chain complex for the Champanerkar-Kofman spanning tree model of reduced Khovanov homology.
\end{abstract}

\date{March 19, 2012.  Published in the Journal of Knot Theory and its Ramifications vol. 21, no. 6 (2012), DOI: 10.1142/S0218216512500629}



\thanks{The author was supported by LSU VIGRE grant DMS-0739382 while conducting much of this research.  Initial interest in this project emerged based on work by Abhijit Champanerkar and Ilya Kofman  \cite{ChKo:mu} and \cite{ChKo:sp}.  The author had several helpful discussions on this subject with Oliver Dasbach and Neal Stoltzfus, and the referee's suggestion to reorganize this work was particularly valuable.}

\maketitle

\section{Introduction}
The computational complexity of the Jones polynomial of an arbitrary knot is $\#P$-hard due to a result by Jaeger, Vertigan, and Welsh \cite{VertWel}.  Loebl and Moffatt \cite{LoeMof} give a \emph{permanent}, or unsigned determinant, formula for the Jones polynomial of any $n$-crossing knot using a $7n\times 7n$ matrix, relying on Monte-Carlo algorithms to estimate it.

However, the Main Theorem \ref{DetJones} of the present paper provides a class of knots for which the Jones polynomial can be obtained by a polynomial-time algorithm given by the determinant of an $n\times n$ matrix.  Proposition \ref{prop:subdividedouble} extends this class beyond pretzel knots, but as is discussed in Remark \ref{remark:Montesinos}, the method below cannot be automatically applied to all Montesinos knots.

This determinant arises from the perfect matching or \emph{dimer} model of a particular bipartite graph $\Gamma$ that can be obtained from a knot diagram.  Futhermore, the unweighted underlying graph $\Gamma$ used here is the same as the one in the dimer model for the Alexander polynomial constructed in another work by the author with Dasbach and Russell \cite{CoDaRu}.

Kasteleyn \cite{Kast} and Temperley and Fisher \cite{TempFish} developed dimers in the 1960's as a tool for studying statistical physics.  For a more formal treatment, see Kenyon's lecture notes \cite{Ken}.  It is because $\Gamma$ is a plane graph that a \emph{Kasteleyn weighting} can be employed to produce a determinant and not just a permanent expansion.

The Jones polynomial is the graded euler characteristic of a complex of bigraded abelian groups called Khovanov homology \cite{Kh:coh}.  Reduced Khovanov homology was described by a state sum of spanning trees by Champanerkar and Kofman \cite{ChKo:mu} and \cite{ChKo:sp} and independently by Wehrli \cite{We:sp}.  The present paper follows some techniques of the first authors; a discussion of how the results in this paper might be useful in this context can be found in Section \ref{sec:Corollaries} below.

Results of a different nature on the same subject of the Jones polynomial of pretzel knots can be found in work by Landvoy \cite{Land} and Jin and Zhang \cite{JinZha}.  Inverse q-determinant formulas for the colored Jones polynomial of a knot are given in work by Huynh and Le \cite{HuynLe} and Garoufalidis and Loebl \cite{GarLoe} independently.





The algorithmic construction in the present paper of the $n\times n$ matrix for a given pretzel knot $K=P(n_1,\ldots,n_k)$ with $n=|n_1|+\ldots+|n_k|$ crossings is as follows.

For each $|n_i|\times(|n_i|-1)$ block, put $L$ in each entry in the main diagonal and $D$ in each entry in the sub-diagonal, and line up these blocks corner-to-corner.  Add a column with $L$ in the $|n_1|$-st position and $D$ in the first row of each of the remaining blocks.  Follow with $k-1$ columns, where the $i$-th of these colums has nonzero entries filling up the entire $i$-th and ($i+1$)-st blocks with a single $\ell$ followed by $d$ exactly $(|n_i|-1)+|n_{i+1}|$ times.

This is the \textit{unenhanced activity matrix} associated to $K$ with entries weighted by activity $\alpha$, technically for all $n_i>0$.  For negative $n_i$, all entries in the corresponding row blocks must be converted to the barred versions $\overline{L}$, $\overline{D}$, $\overline{\ell}$, and $\overline{d}$.  The \textit{enhanced activity matrix} includes an additional writhe weight $w$ that multiplies each row by a monomial depending on its sign, and the signed versions of these matrices have a Kasteleyn weighting $\kappa$ that intersperses some negative signs.

The Main Theorem \ref{DetJones} states that the determinant of this matrix gives the Jones polynomial of any pretzel knot $K$ after suitable evaluations of $L$, $D$, $\ell$, $d$, $\overline{L}$, $\overline{D}$, $\overline{\ell}$, and $\overline{d}$, which are Tutte's activity letters (see for example Chapter X in Bollob{\'a}s \cite{Bol}).  The evaluations are given in Table \ref{tab:ActivityEvaluations} of Section \ref{sec:Background}.


Before evaluation, the terms of the determinant expansion give the complete list of activity words of the spanning trees of the (signed) Tait graph $G$ associated with the knot diagram $D$.  These are verified using some technical lemmas in Section \ref{sec:lemmas}.  It should be noted that although the proof of Lemma \ref{LdD} is similar to the proof of the main theorem in work by Diao, Hetyei, and Hinson \cite{DHH}, this work was done independently of theirs.

Overlay $G$ and its plane dual $G^*$ to obtain a bipartite graph $\widehat{\Gamma}$; delete two vertices to obtain a graph $\Gamma$ with equal-sized vertex sets.  There is a bijection between rooted spanning trees $S$ of $G$ and perfect matchings $\mu$ in $\Gamma$.  This is explored in detail in Section \ref{sec:Construction}.




Three examples appear in Section \ref{sec:Examples}:  the simplest nontrivial knot, the trefoil, given as the ($1,1,1$)-pretzel knot; the first non-alternating knot $8_{19}$ appearing on the Rolfsen Knot Table \cite{knotatlas} given as the ($-2,3,3$)-pretzel knot; and the ($-2,3,7$)-pretzel knot that appears in 3-manifold constructions, e.g. by Jun \cite{Jun}.

\section{The Jones polynomial and Tutte's notion of activity}
\label{sec:Background}

\begin{figure}[h]
\begin{center}
\begin{tikzpicture}
\draw[loosely dashed] (.5,.5) circle (.7cm);
\draw[-] (0,0) -- (1,1);
\draw (1,0) -- (.55,.45);
\draw[-] (0,1) -- (.45,.55);

\draw[loosely dashed] (2.5,.5) circle (.7cm);
\draw[-] (2,0) arc (-45:45:.7cm);
\draw[-] (3,1) arc (135:225:.7cm);

\draw[loosely dashed] (4.5,.5) circle (.7cm);
\draw[-] (5,0) arc (45:135:.7cm);
\draw[-] (4,1) arc (225:315:.7cm);

\end{tikzpicture}
	\caption{Local differences between $L$, $L_0$, and $L_\infty$, respectively.}
	\label{fig:CrossingSmoothings}
\end{center}
\end{figure}
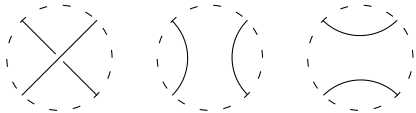

Given an unoriented crossing depicted locally in a link diagram $L$, let $L_0$ and $L_\infty$ be the two smoothings, also called the $A$- and the $B$-smoothings, as in Fig. \ref{fig:CrossingSmoothings}.



Recall the \textit{Kauffman bracket polynomial $\langle L\rangle$} \cite{KauffmanBracket} of a link $L$ defined by
\begin{enumerate}
	\item Smoothing relation:  $\langle L \rangle =A \langle L_0 \rangle + A^{-1} \langle L_\infty \rangle$
	\item Stabilization:  $\langle U\sqcup L\rangle=(-A^2-A^{-2}) \langle L\rangle$
	\item Normalization:  $\langle U\rangle=1$
\end{enumerate}
where $U$ is the unknot and $\sqcup$ is the disjoint union.

\begin{figure}[h]
\begin{center}
\begin{tikzpicture}
\draw[loosely dashed] (.5,.5) circle (.7cm);
\draw[->] (0,0) -- (1,1);
\draw (1,0) -- (.55,.45);
\draw[<-] (0,1) -- (.45,.55);

\draw (-1.5,.5) node {positive};
\draw (4.5,.5) node {negative};

\draw[loosely dashed] (2.5,.5) circle (.7cm);
\draw[<-] (2,1) -- (3,0);
\draw (2,0) -- (2.45,.45);
\draw[<-] (3,1) -- (2.55,.55);

\end{tikzpicture}
	\caption{Positive and negative crossings contribute $+1$ and $-1$ to the writhe.}
	\label{fig:WrithePositiveNegative}
\end{center}
\end{figure}
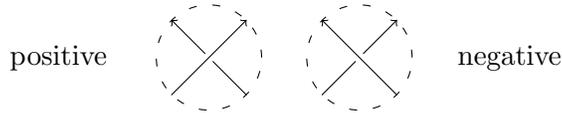

The \textit{writhe} $w(D)$ of an oriented diagram is the sum over all crossings of the evaluation $+1$ for positive crossings and $-1$ for negative crossings as in Fig. \ref{fig:WrithePositiveNegative}.

\begin{definition}
The \textit{Jones polynomial $V_L(t)$} of a link $L$ given an oriented diagram $D$ can be defined via the Kauffman bracket polynomial by
\begin{equation}
V_L(t)=(-A^{-3})^{w(D)}\langle L \rangle,
\end{equation}
where $w(D)$ is the writhe of the 
diagram, along with the substitution $A=t^{-1/4}$.
\end{definition}


The Kauffman bracket polynomial is actually a signed version of the Tutte polynomial of a signed graph $G$ obtained from a link diagram, and the Skein relation (1) comes from the deletion-contraction formula for the Tutte polynomial  (see \cite{Kauff:sign}).

\begin{figure}[h]
\begin{center}

\begin{tikzpicture}
	\fill[gray!50!white] (0,0) -- +(.5,.5) -- +(1,0) -- cycle;
	\fill[gray!50!white] (0,1) -- +(.5,-.5) -- +(1,0) -- cycle;
	\draw[-] (1,0) -- +(-1,1);
	\fill[color=white] (.4,.4) rectangle +(.2,.2);
	\draw[-] (0,0) -- +(1,1);

	\fill[gray!50!white] (2,0) -- +(.5,.5) -- +(1,0) -- cycle;
	\fill[gray!50!white] (2,1) -- +(.5,-.5) -- +(1,0) -- cycle;
	\draw[-] (2,0) -- +(1,1);
	\fill[color=white] (2.4,.4) rectangle +(.2,.2);
	\draw[-] (3,0) -- +(-1,1);

	\draw (-1.5,.5) node {positive};
	\draw (4.5,.5) node {negative};
	
\end{tikzpicture}
	\caption{Crossings determine the sign of the edges in the signed Tait graph.}
	\label{fig:SignedTaitGraphSignsTIKZ}
\end{center}
\end{figure}
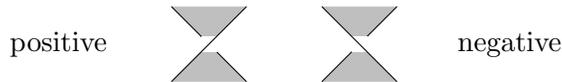

Checkerboard color the regions of a knot diagram black and white.  The \textit{signed Tait graph $G$} associated with a knot diagram has as its vertices the black regions in the checkerboard coloring and as its signed edges the crossings with signs as in Fig.  \ref{fig:SignedTaitGraphSignsTIKZ}.

A spanning tree expansion for the Tutte polynomial can be given in terms of activity words in Tutte's alphabet.  Observe that the activity words themselves are dependent on an ordering of the edges of the graph; invariance of the polynomial is equivalent to proving that any ordering gives the same sum.

\begin{definition}
\label{def:activity}
(Tutte's activity \cite{Bol}) 
To a spanning tree $S$ of a signed graph whose $n$ edges are ordered, Tutte assigns an \emph{activity word} of length $n$ in the alphabet $L,D,\ell,d,\overline{L},\overline{D},\overline{\ell},$ and $\overline{d}$.

A positive edge $e\in S$ is \emph{internally active} (or live, $L$) if it is the lowest numbered edge that reconnects $S-\{e\}$; otherwise it is \emph{internally inactive} (or dead, $D$).

A positive edge $e\notin S$ is \emph{externally active} (or live, $\ell$) if it is the lowest numbered edge in the unique cycle contained in $S\cup \{e\}$; otherwise it is \emph{externally inactive} (or dead, $d$).

Negative edges $\overline{L},\overline{D},\overline{\ell},$ and $\overline{d}$ are defined similarly.

Let $a(e,S)$ be the activity letter associated to the edge $e$ for a specified spanning tree $S$, and let $a(S)$ be the activity word associated to the tree $S$.
\end{definition}

Consider the following small (and very symmetric) example to illustrate activity.

\begin{example}
\label{example:trees}
Let $G$ be the (unsigned, or all positive) graph with two vertices and three parallel edges as in Fig. \ref{fig:trefoiltrees}.  Depicted are the three spanning trees, along with their associated activity words:  ($Ldd)$, $(\ell Dd)$, and $(\ell \ell D)$

\begin{figure}
\begin{center}
\begin{tikzpicture}

\draw[-, rounded corners] (-2.5,2.75) -- (-2.5,3.25) -- (-.5,3.25) -- (-.5,2.75) -- cycle;
\draw[-] (-1.5,3.25) -- (-1.5,2.75);
	\fill[color=black] (-1.5,3.25) circle (3pt);
	\fill[color=black] (-1.5,2.75) circle (3pt);

\draw (-2.8,3) node {$1$};
\draw (-1.8,3) node {$2$};
\draw (-.8,3) node {$3$};

\draw[lightgray] (-.25,2.5) -- (-.25,3.5);

\draw[-, rounded corners, dotted] (3.5,2.75) -- (3.5,3.25) -- (5.5,3.25) -- (5.5,2.75) -- cycle;
\draw[-] (4.5,3.25) -- (4.5,2.75);
	\fill[color=black] (4.5,3.25) circle (3pt);
	\fill[color=black] (4.5,2.75) circle (3pt);

\draw (3.2,3) node {$\ell$};
\draw (4.2,3) node {$D$};
\draw (5.2,3) node {$d$};

\draw[-, rounded corners, dotted] (7.5,2.75) -- (6.5,2.75) -- (6.5,3.25) -- (7.5,3.25);
\draw[-, dotted] (7.5,3.25) -- (7.5,2.75);
\draw[-, rounded corners] (7.5,2.75) -- (8.5,2.75) -- (8.5,3.25) -- (7.5,3.25);
	\fill[color=black] (7.5,3.25) circle (3pt);
	\fill[color=black] (7.5,2.75) circle (3pt);

	\draw (6.2,3) node {$\ell$};
\draw (7.2,3) node {$\ell$};
\draw (8.2,3) node {$D$};

\draw[-, rounded corners] (1.5,2.75) -- (.5,2.75) -- (.5,3.25) -- (1.5,3.25);
\draw[-, dotted] (1.5,3.25) -- (1.5,2.75);
\draw[-, rounded corners, dotted] (1.5,2.75) -- (2.5,2.75) -- (2.5,3.25) -- (1.5,3.25);
	\fill[color=black] (1.5,3.25) circle (3pt);
	\fill[color=black] (1.5,2.75) circle (3pt);

\draw (.2,3) node {$L$};
\draw (1.2,3) node {$d$};
\draw (2.2,3) node {$d$};

\end{tikzpicture}
\end{center}
	\caption{A graph $G$ with its three spanning trees and their associated activity words.}
	\label{fig:trefoiltrees}
\end{figure}
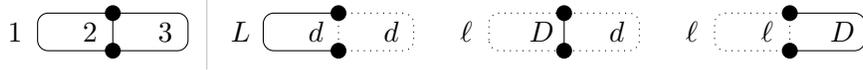
\end{example}


This perspective motivates the following spanning tree expansion for the Jones polynomial, with signed Tait graph $G$ obtained from a checkerboard coloring and signing on edges as in Fig. \ref{fig:SignedTaitGraphSignsTIKZ}.

\begin{theorem}
\label{spJones} (Thistlethwaite \cite{Th}) The Jones polynomial $V_K(t)$ of a knot $K$ with diagram $D$ and signed Tait graph $G$ can be written using the activity evaluations $a(e,S)|_V$ given in Table \ref{tab:ActivityEvaluations} via a spanning tree expansion:
\begin{equation}
V_K(t)=(-A^{-3})^{w(D)}\sum_{S}\prod_{e\in E(G)}a(e,S)|_V,
\end{equation}
summing over all spanning trees $S$ of $G$ and taking the product of the activity of each edge $e$ in $G$.
\end{theorem}

\begin{table}
	\caption{Polynomial activity evaluations.}
{\begin{tabular}{|c|c|c|c|c|c|c|c|c|c|}
\hline
 & & & & & & & & & \\
$a(e,S)$ & activity letter of $e$ w.r.t. $S$ &	$L$ & $D$ & $\ell$ & $d$ & $\overline{L}$ & $\overline{D}$ & $\overline{\ell}$ & $\overline{d}$ \\
 & & & & & & & & & \\
\hline
 & & & & & & & & & \\
 & evaluation for the & & & & & & & & \\
$a(e,S)|_V$ & Jones polynomial &	$-A^{-3}$ & $A$ & $-A^3$ & $A^{-1}$ & $-A^3$ & $A^{-1}$ & $-A^{-3}$ & $A$ \\
 & $V_K(t)$ with $A=t^{-1/4}$ & & & & & & & & \\
 & & & & & & & & & \\
\hline
\end{tabular}}
	\label{tab:ActivityEvaluations}
\end{table}

\section{Technical lemmas classifying activity}
\label{sec:lemmas}

The set of edges between the same two vertices will be referred to as a \textit{parallel edge class}.  The dual of this parallel edge class will be non-standardly referred to as a \textit{2-path}, a connected subgraph $H$ of $G$ where all but two vertices have valence two in both $H$ and $G$, and where these two vertices have valence one in $H$ and any valence in $G$.  Another way to view a 2-path is by the repeated subdivision of an edge.

The following lemmas will be used to determine the activity words of spanning trees on 2-paths and parallel edge classes in signed graphs whose edges are ordered.  

\begin{lemma}
\label{LdD}
(Classification of activity on 2-paths) 
Suppose that a 2-path $P_{k+1}$ of $k+1$ positive edges indexed sequentially by $e_i,\ldots,e_{i+k}$ for some $k>0$ belongs to a graph $G$ such that all interior vertices of the path have degree exactly two in $G$.  Then when determining the activity word for $G$ given a spanning tree $S$, the portion of the activity word associated with the path $P_{k+1}$ must be one of the following:
\begin{enumerate}
	\item ($L\ldots L$) or ($D\ldots D$) when including all edges of $P_{k+1}$ in $S$, or
	\item ($L\ldots LdD\ldots D$) when omitting only the edge $e_{i+j}$ for $0< j\leq k$, or
	\item ($\ell D\ldots D$) or ($d D\ldots D$) when omitting only the first edge.
\end{enumerate}
\end{lemma}

The following proof is similar to the one given for the main theorem in work by Diao, Hetyei, and Hinson \cite{DHH}, although this work was done independently of theirs.

\begin{proof}
First suppose that the 2-path $P_{k+1}$ is contained in the spanning tree $S$; then each of its edges are labelled by $L$ or $D$.  Supposing the edge $e_i$ is labelled $L$, this is the lowest-ordered edge to reconnect the severed tree $S-\{e_i\}$, and so no edge indexed less than $i+1$ can reconnect the severed tree $S-\{e_{i+1}\}$.  Iterate $j$ times to get the 2-path labelled by ($L\ldots L$).  Supposing the edge $e_i$ is labelled $D$, there is an edge indexed less than $i$ that reconnects the severed tree $S-\{e_i\}$, and so this edge also reconnects the severed tree $S-\{e_{i+1}\}$.  Iterate $j$ times to get the 2-path labelled by ($D\ldots D$).  Take note that only these two possibilities arise here, as this fact is used below.

Now suppose that the edge $e_{i+j}$ of the 2-path $P_{k+1}$ is not contained in the spanning tree $S$; then it must be labelled by either $\ell$ or $d$.  Unless $j=0$, the cycle contained in $S\cup \{e_{i+j}\}$ must also contain the edge $e_{i+j-1}$, so $e_{i+j}$ must be labelled by $d$.  Note that no	 other edge may be omitted without disconnecting the tree, so the rest of the edges of $P_{k+1}$ must be labelled by $L$ or $D$.  Then by the argument above, both paths $e_i,\ldots,e_{i+j-1}$ and $e_{i+j+1},\ldots,e_{i+k}$ must be one of ($L\ldots L$) or ($D\ldots D$).

The string before the omitted edge $e_{i+j}$ must contain an $L$ (and therefore be ($L\ldots L$) by the argument above) because only the edges $e_{i+j-1}$ and $e_{i+j}$ reconnect the severed tree $S-\{e_{i+j-1}\}$.  The string following the omitted edge $e_{i+j}$ must contain a $D$ (and therefore be ($D\ldots D$) by the argument above) because only the edges $e_{i+j}$ and $e_{i+j+1}$ reconnect the severed tree $S-\{e_{i+j+1}\}$.

When $j=0$ it is possible for $e_i$ to be labelled by either $\ell$ or $d$, but the string following this omitted edge must be ($D\ldots D$) as above.
\end{proof}

The negative edge case is similar, with each activity letter replaced by its barred counterpart.  Furthermore this works for arbitrary signed edges along a 2-path.  The dual version and its negative edge case work, as well, and the proofs are similar.

The next lemma follows from the previous one by duality.

\begin{lemma}
\label{ellDd}
(Classification of activity on parallel edge classes)  
Suppose a parallel edge class $D_{k+1}$ of $k+1$ positive edges indexed sequentially by $e_i,\ldots,e_{i+k}$ for some $k>0$ belongs to a graph $G$ such that no other edges belong to this class.  Then when determining the activity word for $G$ given a spanning tree $S$, the portion of the activity word associated with the parallel edge class $D_{k+1}$ must be one of the following:
\begin{enumerate}
	\item ($\ell\ldots \ell$) or ($d\ldots d$) when omitting all edges of $D_{k+1}$ in $S$, or
	\item ($\ell\ldots \ell Dd\ldots d$) when including only the edge $e_{i+j}$ for $0< j\leq k$, or
	\item ($L d\ldots d$) or ($D d\ldots d$) when including only the first edge.
\end{enumerate}
\end{lemma}


\section{Constructions of graphs and matrices from knots}
\label{sec:Construction}

\subsection{The signed Tait graph $G$ and balanced overlaid Tait graph $\Gamma$}

Recall the signed Tait graph $G$ associated with the checkerboard coloring of the regions of a knot diagram, and that its vertices are the black regions in the coloring and as its signed edges are the crossings with signs as in Fig. \ref{fig:SignedTaitGraphSignsTIKZ}.

Note that there are actually two graphs here: $G$ for the black regions and its dual $G^*$ in the plane for the white regions, where dual edges change signs.

The \textit{overlaid Tait graph $\widehat{\Gamma}$} is the signed bipartite graph whose first vertex set is the set of intersection points of the edges of both Tait graphs and whose second vertex set is the union of the vertex sets of both Tait graphs.  That is, $V(\widehat{\Gamma})=[E(G)\cap E(G^*)]\sqcup [V(G)\sqcup V(G^*)]$.  The edges of this graph are the half-edges of both Tait graphs.  A similar notion is found in work by Huggett and Virdee \cite{HugVir}.

\begin{figure}
\begin{center}
\begin{tikzpicture}

\draw (-1,12.5) node {$(A)$};
\draw (-1,7.5) node {$(B)$};
\draw (-1,2.5) node {$(C)$};


\draw[-, rounded corners, lightgray] [shift={(0,10)}] (5,.5) -- (5,0) -- (0,0) -- (0,.5);

	\draw[-, lightgray] [shift={(0,10)}] (0,.5) -- +(1,1);
	\fill[color=white] [shift={(0,10)}] (.35,.85) rectangle +(.3,.3);
	\draw[-, lightgray] [shift={(0,10)}] (0,1.5) -- +(1,-1);

\draw[-, lightgray] [shift={(0,10)}] (0,1.5) -- (0,2.5);
\draw[-, lightgray] [shift={(0,10)}] (1,1.5) -- (1,2.5);
	
	\draw[-, lightgray] [shift={(0,10)}] (0,2.5) -- +(1,1);
	\fill[color=white] [shift={(0,10)}] (.35,2.85) rectangle +(.3,.3);
	\draw[-, lightgray] [shift={(0,10)}] (0,3.5) -- +(1,-1);

\draw[-, lightgray] [shift={(0,10)}] (1,.5) -- (2,.5);
\draw[-, lightgray] [shift={(0,10)}] (1,3.5) -- (2,3.5);

	\draw[-, lightgray] [shift={(0,10)}] (2,1.5) -- +(1,-1);
	\fill[color=white] [shift={(0,10)}] (2.35,.85) rectangle +(.3,.3);
	\draw[-, lightgray] [shift={(0,10)}] (2,.5) -- +(1,1);

	\draw[-, lightgray] [shift={(0,10)}] (2,2.5) -- +(1,-1);
	\fill[color=white] [shift={(0,10)}] (2.35,1.85) rectangle +(.3,.3);
	\draw[-, lightgray] [shift={(0,10)}] (2,1.5) -- +(1,1);
	
	\draw[-, lightgray] [shift={(0,10)}] (2,3.5) -- +(1,-1);
	\fill[color=white] [shift={(0,10)}] (2.35,2.85) rectangle +(.3,.3);
	\draw[-, lightgray] [shift={(0,10)}] (2,2.5) -- +(1,1);

\draw[-, lightgray] [shift={(0,10)}] (3,.5) -- (4,.5);
\draw[-, lightgray] [shift={(0,10)}] (3,3.5) -- (4,3.5);

	\draw[-, lightgray] [shift={(0,10)}] (4,1.5) -- +(1,-1);
	\fill[color=white] [shift={(0,10)}] (4.35,.85) rectangle +(.3,.3);
	\draw[-, lightgray] [shift={(0,10)}] (4,.5) -- +(1,1);

	\draw[-, lightgray] [shift={(0,10)}] (4,2.5) -- +(1,-1);
	\fill[color=white] [shift={(0,10)}] (4.35,1.85) rectangle +(.3,.3);
	\draw[-, lightgray] [shift={(0,10)}] (4,1.5) -- +(1,1);
	
	\draw[-, lightgray] [shift={(0,10)}] (4,3.5) -- +(1,-1);
	\fill[color=white] [shift={(0,10)}] (4.35,2.85) rectangle +(.3,.3);
	\draw[-, lightgray] [shift={(0,10)}] (4,2.5) -- +(1,1);

\draw[-, rounded corners, lightgray] [shift={(0,10)}] (5,3.5) -- (5,4) -- (0,4) -- (0,3.5);

\draw[-, rounded corners, double] [shift={(0,10)}] (2.5,.25) -- (.5, .25) -- (.5,3.75) -- (2.5,3.75);
\draw[-, rounded corners] [shift={(0,10)}] (2.5,3.75) -- (4.5,3.75) -- (4.5,.25) -- (2.5,.25);

	\fill[color=black] [shift={(0,10)}] (.5,2) circle (3pt);

	\fill[color=black] [shift={(0,10)}] (2.5,3.75) circle (3pt);

	\fill[color=black] [shift={(0,10)}] (2.5,2.5) circle (3pt);
	\fill[color=black] [shift={(0,10)}] (2.5,1.5) circle (3pt);
	\fill[color=black] [shift={(0,10)}] (2.5,.25) circle (3pt);
\draw[-] [shift={(0,10)}] (2.5,3.75) -- (2.5,.25);

	\fill[color=black] [shift={(0,10)}] (4.5,2.5) circle (3pt);
	\fill[color=black] [shift={(0,10)}] (4.5,1.5) circle (3pt);


\draw[-, rounded corners, lightgray] [shift={(0,5)}] (5,.5) -- (5,0) -- (0,0) -- (0,.5);

	\draw[-, lightgray] [shift={(0,5)}] (0,.5) -- +(1,1);
	\fill[color=white] [shift={(0,5)}] (.35,.85) rectangle +(.3,.3);
	\draw[-, lightgray] [shift={(0,5)}] (0,1.5) -- +(1,-1);

\draw[-, lightgray] [shift={(0,5)}] (0,1.5) -- (0,2.5);
\draw[-, lightgray] [shift={(0,5)}] (1,1.5) -- (1,2.5);
	
	\draw[-, lightgray] [shift={(0,5)}] (0,2.5) -- +(1,1);
	\fill[color=white] [shift={(0,5)}] (.35,2.85) rectangle +(.3,.3);
	\draw[-, lightgray] [shift={(0,5)}] (0,3.5) -- +(1,-1);

\draw[-, lightgray] [shift={(0,5)}] (1,.5) -- (2,.5);
\draw[-, lightgray] [shift={(0,5)}] (1,3.5) -- (2,3.5);

	\draw[-, lightgray] [shift={(0,5)}] (2,1.5) -- +(1,-1);
	\fill[color=white] [shift={(0,5)}] (2.35,.85) rectangle +(.3,.3);
	\draw[-, lightgray] [shift={(0,5)}] (2,.5) -- +(1,1);

	\draw[-, lightgray] [shift={(0,5)}] (2,2.5) -- +(1,-1);
	\fill[color=white] [shift={(0,5)}] (2.35,1.85) rectangle +(.3,.3);
	\draw[-, lightgray] [shift={(0,5)}] (2,1.5) -- +(1,1);
	
	\draw[-, lightgray] [shift={(0,5)}] (2,3.5) -- +(1,-1);
	\fill[color=white] [shift={(0,5)}] (2.35,2.85) rectangle +(.3,.3);
	\draw[-, lightgray] [shift={(0,5)}] (2,2.5) -- +(1,1);

\draw[-, lightgray] [shift={(0,5)}] (3,.5) -- (4,.5);
\draw[-, lightgray] [shift={(0,5)}] (3,3.5) -- (4,3.5);

	\draw[-, lightgray] [shift={(0,5)}] (4,1.5) -- +(1,-1);
	\fill[color=white] [shift={(0,5)}] (4.35,.85) rectangle +(.3,.3);
	\draw[-, lightgray] [shift={(0,5)}] (4,.5) -- +(1,1);

	\draw[-, lightgray] [shift={(0,5)}] (4,2.5) -- +(1,-1);
	\fill[color=white] [shift={(0,5)}] (4.35,1.85) rectangle +(.3,.3);
	\draw[-, lightgray] [shift={(0,5)}] (4,1.5) -- +(1,1);
	
	\draw[-, lightgray] [shift={(0,5)}] (4,3.5) -- +(1,-1);
	\fill[color=white] [shift={(0,5)}] (4.35,2.85) rectangle +(.3,.3);
	\draw[-, lightgray] [shift={(0,5)}] (4,2.5) -- +(1,1);

\draw[-, rounded corners, lightgray] [shift={(0,5)}] (5,3.5) -- (5,4) -- (0,4) -- (0,3.5);

\draw[-, double] [shift={(0,5)}] (1.5,2) -- (4.5,2);
\draw[-, rounded corners] [shift={(0,5)}](.5,3) -- (1,3) -- (1.5,2);
\draw[-, rounded corners] [shift={(0,5)}](1.5,2) -- (1,1) -- (.5,1);
\draw[-, rounded corners, double] [shift={(0,5)}](2.5,3) -- (2,3) -- (1.5,2);
\draw[-, rounded corners, double] [shift={(0,5)}](1.5,2) -- (2,1) -- (2.5,1);
\draw[-, rounded corners, double] [shift={(0,5)}](2.5,3) -- (3,3) -- (3.5,2);
\draw[-, rounded corners, double] [shift={(0,5)}](3.5,2) -- (3,1) -- (2.5,1);
\draw[-, rounded corners, double] [shift={(0,5)}](4.5,3) -- (4,3) -- (3.5,2);
\draw[-, rounded corners, double] [shift={(0,5)}](3.5,2) -- (4,1) -- (4.5,1);

\draw[-, rounded corners] [shift={(0,5)}](.5,3) -- (-.25,3) -- (-.25,4.25) -- (2.5,4.25);
\draw[-, rounded corners, double] [shift={(0,5)}](4.5,3) -- (5.25,3) -- (5.25,4.25) -- (2.5,4.25);

\draw[-, rounded corners, double] [shift={(0,5)}](4.5,2) -- (5.5,2) -- (5.5,4.5) -- (3,4.5) -- (2.5,4.25);

\draw[-, rounded corners] [shift={(0,5)}](.5,1) -- (-.5,1) -- (-.5,4.5) -- (2,4.5) -- (2.5,4.25);
\draw[-, rounded corners, double] [shift={(0,5)}](4.5,1) -- (5.75,1) -- (5.75,4.75) -- (3,4.75) -- (2.5,4.25);

	\fill[color=white] [shift={(0,5)}](1.5,2) circle (3pt);
	\draw [shift={(0,5)}](1.5,2) circle (3pt);
	\fill[color=white] [shift={(0,5)}](3.5,2) circle (3pt);
	\draw [shift={(0,5)}](3.5,2) circle (3pt);
	\fill[color=white] [shift={(0,5)}](2.5,4.25) circle (3pt);
	\draw [shift={(0,5)}](2.5,4.25) circle (3pt);



\draw[-, rounded corners, lightgray] (5,.5) -- (5,0) -- (0,0) -- (0,.5);

	\draw[-, lightgray] (0,.5) -- +(1,1);
	\fill[color=white] (.35,.85) rectangle +(.3,.3);
	\draw[-, lightgray] (0,1.5) -- +(1,-1);

\draw[-, lightgray] (0,1.5) -- (0,2.5);
\draw[-, lightgray] (1,1.5) -- (1,2.5);
	
	\draw[-, lightgray] (0,2.5) -- +(1,1);
	\fill[color=white] (.35,2.85) rectangle +(.3,.3);
	\draw[-, lightgray] (0,3.5) -- +(1,-1);

\draw[-, lightgray] (1,.5) -- (2,.5);
\draw[-, lightgray] (1,3.5) -- (2,3.5);

	\draw[-, lightgray] (2,1.5) -- +(1,-1);
	\fill[color=white] (2.35,.85) rectangle +(.3,.3);
	\draw[-, lightgray] (2,.5) -- +(1,1);

	\draw[-, lightgray] (2,2.5) -- +(1,-1);
	\fill[color=white] (2.35,1.85) rectangle +(.3,.3);
	\draw[-, lightgray] (2,1.5) -- +(1,1);
	
	\draw[-, lightgray] (2,3.5) -- +(1,-1);
	\fill[color=white] (2.35,2.85) rectangle +(.3,.3);
	\draw[-, lightgray] (2,2.5) -- +(1,1);

\draw[-, lightgray] (3,.5) -- (4,.5);
\draw[-, lightgray] (3,3.5) -- (4,3.5);

	\draw[-, lightgray] (4,1.5) -- +(1,-1);
	\fill[color=white] (4.35,.85) rectangle +(.3,.3);
	\draw[-, lightgray] (4,.5) -- +(1,1);

	\draw[-, lightgray] (4,2.5) -- +(1,-1);
	\fill[color=white] (4.35,1.85) rectangle +(.3,.3);
	\draw[-, lightgray] (4,1.5) -- +(1,1);
	
	\draw[-, lightgray] (4,3.5) -- +(1,-1);
	\fill[color=white] (4.35,2.85) rectangle +(.3,.3);
	\draw[-, lightgray] (4,2.5) -- +(1,1);

\draw[-, rounded corners, lightgray] (5,3.5) -- (5,4) -- (0,4) -- (0,3.5);


\draw[rounded corners] (1.5,2) -- (1,1) -- (.5,1);

\draw[] (2.5,1.5) -- (2.5,1);
\draw[] (2.5,2.5) -- (2.5,2);

\draw[] (4.5,3) -- (4.5,2.5);
\draw[] (4.5,2) -- (4.5,1.5);
\draw[rounded corners] (4.5,1) -- (4.5,.25) -- (2.5,.25);


\draw[-, rounded corners] (.5,3) -- (-.25,3) -- (-.25,4.25) -- (2.5,4.25);
\draw[-, rounded corners] (4.5,3) -- (5.25,3) -- (5.25,4.25) -- (2.5,4.25);

\draw[-, rounded corners] (4.5,2) -- (5.5,2) -- (5.5,4.5) -- (3,4.5) -- (2.5,4.25);

\draw[-, rounded corners] (.5,1) -- (-.5,1) -- (-.5,4.5) -- (2,4.5) -- (2.5,4.25);
\draw[-, rounded corners] (4.5,1) -- (5.75,1) -- (5.75,4.75) -- (3,4.75) -- (2.5,4.25);


\draw[-, rounded corners] (.5,1) --(.5,.25) -- (2.5,.25);
\draw[-] (.5,3) -- (.5,1);
	\fill[color=white] (.5,2) circle (3pt);
	\draw (.5,2) circle (3pt);

\draw[-] (2.5,3) -- (2.5,2.5);
\draw[-] (2.5,2) -- (2.5,1.5);
\draw[-] (2.5,1) -- (2.5,.25);

\draw[-, rounded corners] (.5,3) -- (.5,3.75) -- (4.5,3.75) -- (4.5,3);
\draw[-] (2.5,3) -- (2.5,3.75);

	\fill[color=white] (2.5,4.25) circle (3pt);
	\draw (2.5,4.25) circle (3pt);
	\fill[color=white] (2.5,3.75) circle (3pt);
	\draw (2.5,3.75) circle (3pt);

	\fill[color=white] (2.5,2.5) circle (3pt);
	\draw (2.5,2.5) circle (3pt);
	\fill[color=white] (2.5,1.5) circle (3pt);
	\draw (2.5,1.5) circle (3pt);
	\fill[color=white] (2.5,.25) circle (3pt);
	\draw (2.5,.25) circle (3pt);

\draw[-] (4.5,2.5) -- (4.5,2);
\draw[-] (4.5,1.5) -- (4.5,1);
	\fill[color=white] (4.5,2.5) circle (3pt);
	\draw (4.5,2.5) circle (3pt);
	\fill[color=white] (4.5,1.5) circle (3pt);
	\draw (4.5,1.5) circle (3pt);


	\fill[color=black] (.5,3) circle (3pt);
	\fill[color=black] (.5,1) circle (3pt);

	\fill[color=black] (2.5,3) circle (3pt);
	\fill[color=black] (2.5,2) circle (3pt);
	\fill[color=black] (2.5,1) circle (3pt);

	\fill[color=black] (4.5,3) circle (3pt);
	\fill[color=black] (4.5,2) circle (3pt);
	\fill[color=black] (4.5,1) circle (3pt);


\draw[-] (1.5,2) -- (4.5,2);
\draw[-, rounded corners] (.5,3) -- (1,3) -- (1.5,2);
\draw[-, rounded corners] (2.5,3) -- (2,3) -- (1.5,2);
\draw[-, rounded corners] (1.5,2) -- (2,1) -- (2.5,1);
\draw[-, rounded corners] (2.5,3) -- (3,3) -- (3.5,2);
\draw[-, rounded corners] (3.5,2) -- (3,1) -- (2.5,1);
\draw[-, rounded corners] (4.5,3) -- (4,3) -- (3.5,2);
\draw[-, rounded corners] (3.5,2) -- (4,1) -- (4.5,1);

	\fill[color=white] (1.5,2) circle (3pt);
	\draw (1.5,2) circle (3pt);
	\fill[color=white] (3.5,2) circle (3pt);
	\draw (3.5,2) circle (3pt);

\end{tikzpicture}
	\caption{(A) The Tait graph $G$ for $8_{19}$, (B) its dual $G^*$, and (C) the overlaid Tait graph $\widehat{\Gamma}$.}
	\label{fig:graphs}
\end{center}
\end{figure}
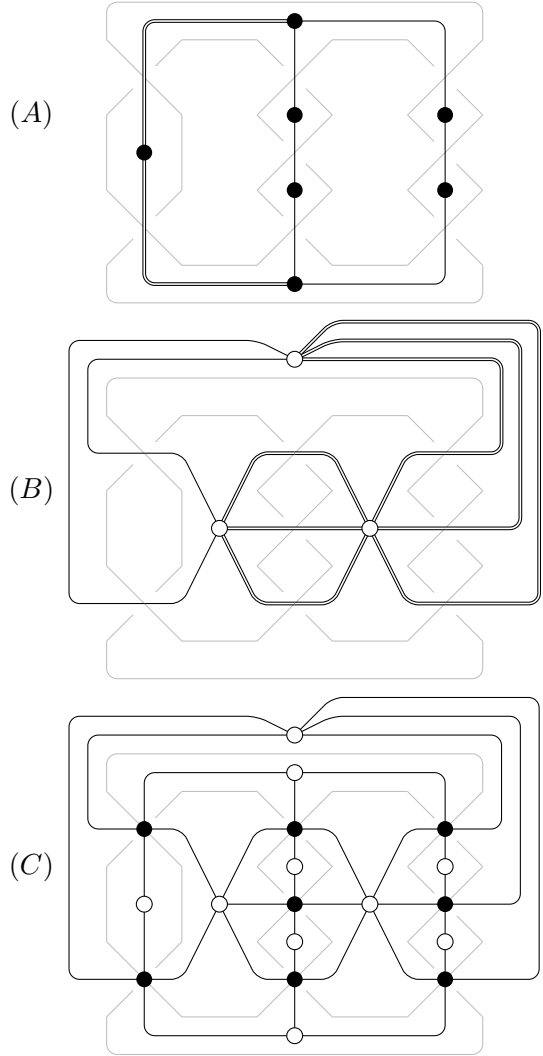

Fig. \ref{fig:graphs} shows (A) the signed Tait graph $G$, (B) it's dual $G^*$, and (C) the overlaid Tait graph $\widehat{\Gamma}$ for the knot $8_{19}$ considered in Example \ref{example2} below.  The thickened edges in (A) and (B) are negatively signed.

Note that the signs of the edges of the overlaid Tait graph $\widehat{\Gamma}$ do not arise from the signs of the original Tait graphs but may be assigned somewhat arbitrarily according to a Kasteleyn weighting, which is Definition \ref{def:KasteleynWeighting} below.

\begin{figure}
\begin{center}
\begin{pspicture}(0,0)(4,2)

\pcline[linewidth=1 pt, linecolor=lightgray]{->}(0,1)(4,1)
\pcline[linewidth=1 pt, linecolor=lightgray]{-}(1,0)(1,2)
\pcline[linewidth=1 pt, linecolor=lightgray]{-}(3,0)(3,2)

\pscircle[linewidth=1pt, linecolor=black, fillstyle=solid, fillcolor=black](1,1){.15}
\pscircle[linewidth=1pt, linecolor=black, fillstyle=solid, fillcolor=black](3,1){.15}

\pcline[linewidth=1 pt]{-}(1,1)(2,0)
\pcline[linewidth=1 pt]{-}(3,1)(2,0)
\pcline[linewidth=1 pt]{-}(1,1)(2,2)
\pcline[linewidth=1 pt]{-}(3,1)(2,2)

\pscircle[linewidth=1pt, linecolor=black, fillstyle=solid](2,0){.15}
\pscircle[linewidth=1pt, linecolor=black, fillstyle=solid](2,2){.15}

\end{pspicture}
	\caption{A square face of the overlaid Tait graph $\widehat{\Gamma}$.}
	\label{fig:SquareFace}
\end{center}
\end{figure}
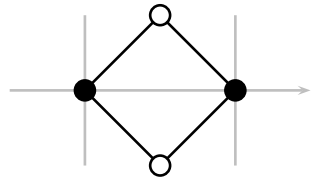

The projection of a knot diagram without its crossing information is called the \textit{projection graph}.  This four-valent graph has as its vertices the crossings of the diagram.  Note that each edge in the projection graph corresponds to a square face in the overlaid Tait graph $\widehat{\Gamma}$, as in Fig. \ref{fig:SquareFace}.

The \textit{balanced overlaid Tait graph $\Gamma$} is obtained from the overlaid Tait graph $\widehat{\Gamma}$ by deleting two vertices from the larger vertex set of $\widehat{\Gamma}$ that lie on the same square face.  This is the graph on which perfect matchings will be considered for the dimer model.  A bipartite graph whose vertex sets are not of the same size, that is, which is \textit{unbalanced}, can have no perfect matchings.

\subsection{The activity matrix $A$}

The incidence matrix of a graph $G$ has its rows labelled by vertices and its columns by edges.  Take the transpose of this matrix for the signed Tait graph $G$ and amalgamate it with the transpose of the incidence matrix for the dual graph $G^*$; thus the columns are partitioned into two sets corresponding to vertices and faces of $G$.  Delete a column of this matrix from each set such that the associated vertex and face are incident.  Call this the \emph{squared incidence matrix} associated with the signed Tait graph $G$.

On the other hand, the adjacency matrix of the balanced overlaid Tait graph $\Gamma$ can be written in block form by
\begin{equation}
\left( \begin{array}{c|c}
0 & M \\
\hline
M^T & 0 \end{array} \right)
\end{equation}
with some submatrix $M$, which shall be called the \textit{bipartite adjacency submatrix} associated with the balanced overlaid Tait graph $\Gamma$.

The matrices are the same as long as the deleted columns of the squared incidence matrix correspond to the vertices of the overlaid Tait graph $\widehat{\Gamma}$ that are deleted to obtain the balanced overlaid Tait graph $\Gamma$.  By construction, the rows of both represent the original crossings of the knot diagram, and the columns of both represent (all but two of) the faces of the projection graph associated to the knot diagram.  For clarity of context, however, both names will be used.

The entries in this matrix are 0's and formal variables in place of the usual 1's.

Recall that in the bipartition of the vertices in the balanced overlaid Tait graph $\Gamma$ the second vertex set is the union of the vertex sets of both the original Tait graph $G$ and its dual $G^*$.  Thus this graph $\Gamma$ is more technically a tripartite graph whose three vertex sets correspond to the (ordered) edges, (all but one of the) vertices, and (all but one of the) faces of the specified Tait graph $G$.  That is, $V(\Gamma)=[E(G)\cap E(G^*)]\sqcup [V(G)] \sqcup [V(G^*)]$.

\begin{definition}
\label{def:activityweighting}
The \textit{activity weighting} $\alpha(\varepsilon)$ on an edge $\varepsilon$ of a balanced overlaid Tait graph $\Gamma$ associated with a knot diagram whose $n$ crossings are ordered is determined by three distinctions:  positive or negative, internal or external, and live or dead.

The activity weighting for an edge incident with a vertex from the first set is positive or negative if the edge corresponding to the crossing in the chosen Tait graph is signed positive or negative, respectively, according to Fig. \ref{fig:SignedTaitGraphSignsTIKZ}.

The activity weighting for an edge is internal or external depending on whether it is incident with a vertex in the second or third set, respectively, according to the sets above.  Note that every edge is incident with a vertex in the first set, so this is indeed a partition.

The activity weighting for an edge is live or dead depending on whether or not it connects the lowest-ordered vertex from the first set to the other vertex with which it is incident.  That is, given a vertex in the second or third set, label all of the edges around it as dead except for the one incident with the lowest-ordered vertex in the first set, and this last edge is labelled live.

These three choices determine the activity letter assigned as the weighting of the edge.
\end{definition}

The evaluations $\alpha(\varepsilon)|_V$ of $\alpha(\varepsilon)$ follow those given in Table \ref{tab:ActivityEvaluations}.  Note the distinction between the activity evaluations $\alpha(\varepsilon)$ of edges $\varepsilon$ in the balanced overlaid Tait graph $\Gamma$ and the activity evaluations $a(e,S)$ of edges $e$ in the Tait graph $G$; the former are not dependent on a choice of spanning tree $S$!

The \textit{writhe weighting} $w(\varepsilon)$ on an edge $\varepsilon$ of a balanced overlaid Tait graph $\Gamma$ associated with a knot diagram is either $(-A)^{-3}$ or $(-A)^{3}$, depending on the positive or negative contribution, respectively, to the writhe by the crossing whose associated vertex the edge $\varepsilon$ is incident with.

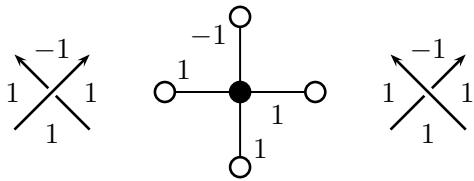
\begin{figure}
\begin{center}
\begin{pspicture}(0,-.5)(6.5,1.5)
\pcline[linewidth=1 pt]{->}(0,0)(1,1)
\pcline[linewidth=1 pt]{<-}(0,1)(.45,.55)
\pcline[linewidth=1 pt]{-}(1,0)(.55,.45)

\uput[0](.75,.5){$1$}
\uput[90](.5,.75){$-1$}
\uput[180](.25,.5){$1$}
\uput[270](.5,.25){$1$}

\pcline[linewidth=.75 pt]{-}(2,.5)(4,.5)
\pcline[linewidth=.75 pt]{-}(3,-.5)(3,1.5)

\uput[270](3.5,.5){$1$}
\uput[180](3,1.25){$-1$}
\uput[90](2.25,.5){$1$}
\uput[0](3,-.25){$1$}

\pscircle[linewidth=1pt, linecolor=black, fillstyle=solid, fillcolor=black](3,.5){.15}

\pscircle[linewidth=1pt, linecolor=black, fillstyle=solid](2,.5){.15}
\pscircle[linewidth=1pt, linecolor=black, fillstyle=solid](4,.5){.15}
\pscircle[linewidth=1pt, linecolor=black, fillstyle=solid](3,-.5){.15}
\pscircle[linewidth=1pt, linecolor=black, fillstyle=solid](3,1.5){.15}

\pcline[linewidth=1 pt]{->}(6,0)(5,1)
\pcline[linewidth=1 pt]{<-}(6,1)(5.55,.55)
\pcline[linewidth=1 pt]{-}(5,0)(5.45,.45)

\uput[0](5.75,.5){$1$}
\uput[90](5.5,.75){$-1$}
\uput[180](5.25,.5){$1$}
\uput[270](5.5,.25){$1$}

\end{pspicture}
	\caption{Kauffman's trick $\kappa(\varepsilon)$ to produce a Kasteleyn weighting.}
	\label{fig:KauffmansTrick}
\end{center}
\end{figure}

Kauffman's signing convention (which appears in \cite{Kauff2}, the republication of \cite{Kauff}) described in Fig. \ref{fig:KauffmansTrick} provides a canonical way to distribute signs to the edges of the balanced overlaid Tait graph $\Gamma$ coming from a knot diagram.  This shall be called \textit{Kauffman's trick} and shall be denoted $\kappa(\varepsilon)$ for an edge $\varepsilon\in E(\Gamma)$.

\begin{definition}
\label{def:activitymatrix}
The \textit{unenhanced activity matrix} $A$ associated with a knot diagram is the squared incidence matrix or bipartite adjacency submatrix whose non-zero entries are given by the activity weighting $\alpha(\varepsilon)$ above.  Together with the writhe weighting $w(\varepsilon)$ this gives the \textit{enhanced activity matrix}.  The \textit{signed unenhanced (or enhanced) activity matrix} has signs coming from a Kasteleyn weighting given by Kauffman's trick $\kappa(\varepsilon)$.

The activity weighting amounts to the following weighting rules.  Ordered rows associated with the original ordered crossings (vertices of the first set) contain only positive or only negative letters following the sign of the original crossing in the specific Tait graph considered.  Columns associated with the vertices of the second set are internal, and columns associated with the vertices of the third set are external.  The first non-zero entry in a column is live; the rest are dead.
\end{definition}

Note that order will be significant only for the rows and not for the columns of this matrix; however, it is advantageous to keep the partition of the columns.


\subsection{Rooted spanning trees of $G$ and perfect matchings of $\Gamma$}

The remainder of this section interprets the determinant expansion of the matrix defined above with formal variables before the activity and writhe weightings.  First consider the expansion from the unsigned version.

The \textit{permanent} or \textit{unsigned determinant} of an $n\times n$ matrix $M=(m_{ij})$ is
\begin{equation}
perm(M)=\sum_{\sigma\in S_n} \prod_{i=1}^n m_{i\sigma(i)},
\end{equation}
summing over all permuations $\sigma$ in the symmetric group $S_n$.

\begin{proposition}
\label{prop:PermDimer}
The non-zero terms in the permanent expansion of a bipartite adjacency submatrix associated with a balanced bipartite graph give the complete list of perfect matchings of the graph.
\end{proposition}

\begin{proof}
Each term in the permanent expansion is a permutation $\sigma$ matching each vertex $i$ in the first vertex set to a vertex $\sigma(i)$ in the second vertex set.
\end{proof}

Thus the important object here is the permanent and not the determinant.  However, for the purpose of easier calculation, the following notion can be used to switch back and forth between the signed and unsigned versions.

\begin{definition}
\label{def:KasteleynWeighting}
A \textit{Kasteleyn weighting} of a plane bipartite graph is a choice of sign for each edge such that the number of negatives around a particular face is
\begin{itemize}
	\item odd if the face has length 0 mod 4 or
	\item even if the face has length 2 mod 4.
\end{itemize}
\end{definition}

Note that this definition only makes sense for a \textit{plane graph}, that is, an abstract graph together with its plane embedding.

\begin{lemma}
\label{lemma:WeightingEdgeDelete}
If an edge is deleted from a graph with a Kasteleyn weighting, the resulting graph still has a Kasteleyn weighting.
\end{lemma}

\begin{proof}
The deletion of an edge incident with two faces of length $f_1$ and $f_2$ results in a new face of length $f_1+f_2-2$.  The number of negatives in this new face changes by 0 or 2 (an even number) compared with the sum of the number of negatives in $f_1$ and $f_2$.

It is left to the reader to check that the four parity cases preserve the Kasteleyn weighting.
\end{proof}


\begin{proposition}
Kauffman's trick $\kappa(\varepsilon)$ provides a Kasteleyn weighting for the balanced overlaid Tait graph $\Gamma$ coming from an oriented knot diagram.
\end{proposition}

\begin{proof}
Recall that each of the faces in the overlaid Tait graph $\widehat{\Gamma}$ is a square as in Fig. \ref{fig:SquareFace}.  The assigning of a negative edge according to Fig. \ref{fig:KauffmansTrick} affects exactly one of the northwest and southwest sides of this square.  Thus exactly one edge of every square face is negatively signed.

By Lemma \ref{lemma:WeightingEdgeDelete}, the edge deletions that result in the balanced overlaid Tait graph $\Gamma$ do not affect this weighting.
\end{proof}

In the discussion below, the signs coming from Kauffman's trick $\kappa(\varepsilon)$ giving a Kasteleyn weighting will not appear on the graph; instead these signs will occur only in the associated terms of the matrix.

\begin{proposition}
\label{prop:DetDimer}
Suppose a balanced bipartite graph is a plane graph whose bipartite adjacency submatrix has entries given signs by a Kasteleyn weighting.  Then the non-zero terms in the determinant expansion of the signed bipartite adjacency submatrix associated with the graph give the complete list of perfect matchings of the unsigned graph, up to an overall sign.
\end{proposition}

\begin{proof}
By Proposition \ref{prop:PermDimer}, only the sign of each term in the determinant expansion needs to be checked against the signs of the edges in the perfect matchings.  It is enough to demonstrate the sign difference between any two terms and the sign difference between the two corresponding perfect matchings are indeed the same.

Suppose two permutations that do not give zero terms differ by exactly one transposition.  This holds if and only if there are four non-zero terms arranged as corners of a rectangle in the matrix.  This holds if an only if there are two vertices from each of the two vertex sets incident with both of the vertices in the other two set, that is, if and only if there is a square face in the graph.

Since the face has an odd number of signs by the Kasteleyn weighting, the two perfect matchings, which differ only on the opposite sides of this square, must have opposite signs.

Any two permutations differ in some number of transpositions, so this can be extended to all terms.
\end{proof}

\begin{proposition}
\label{prop:DimerTrees}
Given a knot diagram and the choice of two omitted faces in the projection graph associated with the knot digram, there is a bijection between perfect matchings of the balanced overlaid Tait graph $\Gamma$ associated with the knot diagram and rooted spanning trees of the Tait graph $G$ associated with the knot diagram.
\end{proposition}

\begin{proof}
By Proposition \ref{prop:DetDimer}, there is a bijection between perfect matchings of the balanced overlaid Tait graph $\Gamma$ associated with the knot diagram and the non-zero terms in the determinant expansion of the bipartite adjacency submatrix.  By construction the bipartite adjacency submatrix is the squared incidence matrix of the Tait graph $G$ associated with the knot diagram.  By Kirchhoff's matrix tree theorem, the non-zero terms of the determinant expansion of the squared incidence matrix give the complete list of rooted spanning trees of the Tait graph $G$ associated with the knot diagram, along with the complementary set of rooted spanning trees in the dual graph $G^*$.
\end{proof}

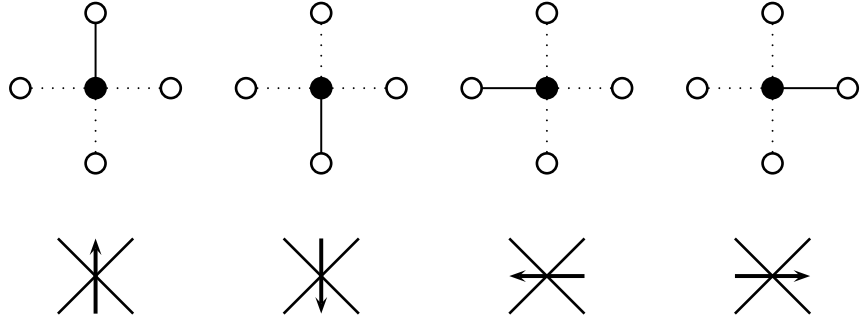
\begin{figure}
\begin{center}
\begin{pspicture}(0,0)(11,4)

\pscircle[linewidth=1pt, linecolor=black, fillstyle=solid, fillcolor=black](1,3){.15}
\pcline[linewidth=.75 pt, linestyle=dotted]{-}(0,3)(1,3)
\pcline[linewidth=.75 pt, linestyle=dotted]{-}(1,3)(2,3)
\pcline[linewidth=.75 pt, linestyle=dotted]{-}(1,2)(1,3)
\pcline[linewidth=.75 pt, linestyle=solid]{-}(1,3)(1,4)
\pscircle[linewidth=1pt, linecolor=black, fillstyle=solid](0,3){.15}
\pscircle[linewidth=1pt, linecolor=black, fillstyle=solid](2,3){.15}
\pscircle[linewidth=1pt, linecolor=black, fillstyle=solid](1,2){.15}
\pscircle[linewidth=1pt, linecolor=black, fillstyle=solid](1,4){.15}

\pscircle[linewidth=1pt, linecolor=black, fillstyle=solid, fillcolor=black](4,3){.15}
\pcline[linewidth=.75 pt, linestyle=dotted]{-}(3,3)(4,3)
\pcline[linewidth=.75 pt, linestyle=dotted]{-}(4,3)(5,3)
\pcline[linewidth=.75 pt, linestyle=solid]{-}(4,2)(4,3)
\pcline[linewidth=.75 pt, linestyle=dotted]{-}(4,3)(4,4)
\pscircle[linewidth=1pt, linecolor=black, fillstyle=solid](3,3){.15}
\pscircle[linewidth=1pt, linecolor=black, fillstyle=solid](5,3){.15}
\pscircle[linewidth=1pt, linecolor=black, fillstyle=solid](4,2){.15}
\pscircle[linewidth=1pt, linecolor=black, fillstyle=solid](4,4){.15}

\pscircle[linewidth=1pt, linecolor=black, fillstyle=solid, fillcolor=black](7,3){.15}
\pcline[linewidth=.75 pt, linestyle=solid]{-}(6,3)(7,3)
\pcline[linewidth=.75 pt, linestyle=dotted]{-}(7,3)(8,3)
\pcline[linewidth=.75 pt, linestyle=dotted]{-}(7,2)(7,3)
\pcline[linewidth=.75 pt, linestyle=dotted]{-}(7,3)(7,4)
\pscircle[linewidth=1pt, linecolor=black, fillstyle=solid](6,3){.15}
\pscircle[linewidth=1pt, linecolor=black, fillstyle=solid](8,3){.15}
\pscircle[linewidth=1pt, linecolor=black, fillstyle=solid](7,2){.15}
\pscircle[linewidth=1pt, linecolor=black, fillstyle=solid](7,4){.15}

\pscircle[linewidth=1pt, linecolor=black, fillstyle=solid, fillcolor=black](10,3){.15}
\pcline[linewidth=.75 pt, linestyle=dotted]{-}(9,3)(10,3)
\pcline[linewidth=.75 pt, linestyle=solid]{-}(10,3)(11,3)
\pcline[linewidth=.75 pt, linestyle=dotted]{-}(10,2)(10,3)
\pcline[linewidth=.75 pt, linestyle=dotted]{-}(10,3)(10,4)
\pscircle[linewidth=1pt, linecolor=black, fillstyle=solid](9,3){.15}
\pscircle[linewidth=1pt, linecolor=black, fillstyle=solid](11,3){.15}
\pscircle[linewidth=1pt, linecolor=black, fillstyle=solid](10,2){.15}
\pscircle[linewidth=1pt, linecolor=black, fillstyle=solid](10,4){.15}

\pcline[linewidth=1 pt]{-}(.5,0)(1.5,1)
\pcline[linewidth=1 pt]{-}(1.5,0)(.5,1)
\pcline[linewidth=1.5 pt]{->}(1,0)(1,1)

\pcline[linewidth=1 pt]{-}(3.5,0)(4.5,1)
\pcline[linewidth=1 pt]{-}(4.5,0)(3.5,1)
\pcline[linewidth=1.5 pt]{<-}(4,0)(4,1)

\pcline[linewidth=1 pt]{-}(6.5,0)(7.5,1)
\pcline[linewidth=1 pt]{-}(7.5,0)(6.5,1)
\pcline[linewidth=1.5 pt]{<-}(6.5,.5)(7.5,.5)

\pcline[linewidth=1 pt]{-}(9.5,0)(10.5,1)
\pcline[linewidth=1 pt]{-}(10.5,0)(9.5,1)
\pcline[linewidth=1.5 pt]{->}(9.5,.5)(10.5,.5)

\end{pspicture}
	\caption{The correspondence between edges $\varepsilon$ in the overlaid Tait graph $\widehat{\Gamma}$ and directed edges $e$ in the (directed) Tait graph $G$.}
	\label{fig:CrossingTaitOverlaidTait}
\end{center}
\end{figure}

To a crossing in the knot diagram, there are exactly four configurations of edges that can be incident with the associated vertex in the overlaid Tait graph $\widehat{\Gamma}$, and there are exactly four configurations of directed edges that can be associated with it in the (directed) Tait graph $G$.  It is easy to see the relationship between these, as depicted in Fig. \ref{fig:CrossingTaitOverlaidTait}.

\section{Main Results}
\label{sec:MainResults}


\subsection{Pretzel knots}

Let $P(n_1,n_2,\ldots,n_k)$ be the \textit{pretzel knot} with $k$ columns, each of $|n_i|\in\mathbb{N}$ crossings, where the sign of $n_i\neq0$ determines the signs of the crossings in the column.

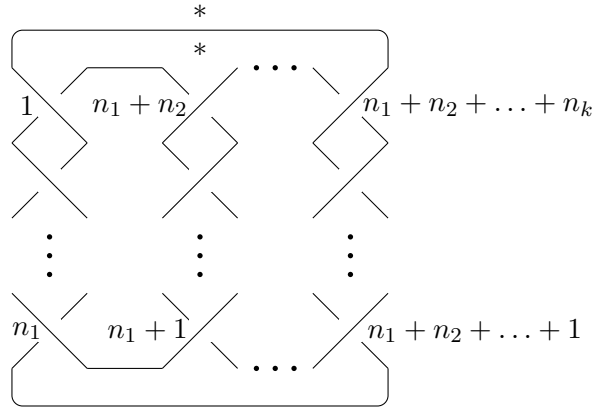
\begin{figure}
\begin{center}
\begin{tikzpicture}

\draw[-, rounded corners]  (5,.5) -- (5,0) -- (0,0) -- (0,.5);

	\draw[-]  (0,.5) -- +(1,1);
	\fill[color=white]  (.35,.85) rectangle +(.3,.3);
	\draw[-]  (0,1.5) -- +(1,-1);
\draw  (.2,1) node {$n_1$};

	\fill[color=black] (.5,2.25) circle (1pt);
	\fill[color=black] (.5,2) circle (1pt);
	\fill[color=black] (.5,1.75) circle (1pt);

	\draw[-]  (0,2.5) -- +(1,1);
	\fill[color=white]  (.35,2.85) rectangle +(.3,.3);
	\draw[-]  (0,3.5) -- +(1,-1);

	\draw[-]  (0,3.5) -- +(1,1);
	\fill[color=white]  (.35,3.85) rectangle +(.3,.3);
	\draw[-]  (0,4.5) -- +(1,-1);
\draw  (.2,4) node {$1$};

\draw[-]  (1,.5) -- (2,.5);
\draw[-]  (1,4.5) -- (2,4.5);

	\draw[-]  (2,1.5) -- +(1,-1);
	\fill[color=white]  (2.35,.85) rectangle +(.3,.3);
	\draw[-]  (2,.5) -- +(1,1);
\draw  (1.8,1) node {$n_1+1$};

	\fill[color=black] (2.5,2.25) circle (1pt);
	\fill[color=black] (2.5,2) circle (1pt);
	\fill[color=black] (2.5,1.75) circle (1pt);
	
	\draw[-]  (2,3.5) -- +(1,-1);
	\fill[color=white]  (2.35,2.85) rectangle +(.3,.3);
	\draw[-]  (2,2.5) -- +(1,1);

	\draw[-]  (2,4.5) -- +(1,-1);
	\fill[color=white]  (2.35,3.85) rectangle +(.3,.3);
	\draw[-]  (2,3.5) -- +(1,1);
\draw  (1.7,4) node {$n_1+n_2$};

	\fill[color=black] (3.25,.5) circle (1pt);
	\fill[color=black] (3.5,.5) circle (1pt);
	\fill[color=black] (3.75,.5) circle (1pt);

	\fill[color=black] (3.25,4.5) circle (1pt);
	\fill[color=black] (3.5,4.5) circle (1pt);
	\fill[color=black] (3.75,4.5) circle (1pt);

	\draw[-]  (4,1.5) -- +(1,-1);
	\fill[color=white]  (4.35,.85) rectangle +(.3,.3);
	\draw[-]  (4,.5) -- +(1,1);
\draw  (6.15,1) node {$n_1+n_2+\ldots+1$};

	\fill[color=black] (4.5,2.25) circle (1pt);
	\fill[color=black] (4.5,2) circle (1pt);
	\fill[color=black] (4.5,1.75) circle (1pt);
	
	\draw[-]  (4,3.5) -- +(1,-1);
	\fill[color=white]  (4.35,2.85) rectangle +(.3,.3);
	\draw[-]  (4,2.5) -- +(1,1);

	\draw[-]  (4,4.5) -- +(1,-1);
	\fill[color=white]  (4.35,3.85) rectangle +(.3,.3);
	\draw[-]  (4,3.5) -- +(1,1);
\draw  (6.2,4) node {$n_1+n_2+\ldots+n_k$};

\draw[-, rounded corners]  (5,4.5) -- (5,5) -- (0,5) -- (0,4.5);

	\draw (2.5,5.25) node {$*$};
	\draw (2.5,4.75) node {$*$};

\end{tikzpicture}
\end{center}
	\caption{The ordering of the crossings and the two omitted regions in pretzel knot $P$.}
	\label{fig:pretzel2TIKZ}
\end{figure}

For the sake of fixing notation, let $G$ be the signed Tait graph and $\Gamma$ the balanced overlaid Tait graph for a diagram $D$ of a pretzel knot $K$ whose $n=|n_1|+\ldots+|n_k|$ crossings are labelled from left to right in columns and downward on the first column and then upward on the remaining columns, with the two omitted faces of the projection graph corresponding to the universal face and the upper deck supported by the columns as in Fig. \ref{fig:pretzel2TIKZ}.

For the diagram of the $(n_1,\ldots,n_k)$-pretzel knot in Fig. \ref{fig:pretzel2TIKZ}, see as an example the following (unsigned) unenhanced activity matrix:
\begin{equation} 
\left(
\begin{array}{cccc|cccc|cccc|c|cccc|c|ccccc}
\overline{L} &   &        &   &   &   &        &   &  &&&&  &   &   &        &   &   & \overline{\ell} &      &&      &      \\
\overline{D} & \overline{L} &        &   &   &   &        &   &  &&&&  &   &   &        &   &   & \overline{d}    &      &&      &      \\
  & \overline{D} & \ddots &   &   &   &        &   &  &&&&  &   &   &        &   &   & \overline{d}    &      &&      &      \\
  &   & \ddots & \overline{L} &   &   &        &   &  &&&&  &   &   &        &   &   & \vdots    &&      &      &      \\
  &   &        & \overline{D} &   &   &        &   &  &&&&  &   &   &        &   & \overline{L} & \overline{d}    &      &&      &      \\
\hline
  &   &        &   & L &   &        &   &  &&&&  &   &   &        &   & D & d    & \ell &&      &      \\
  &   &        &   & D & L &        &   &  &&&&  &   &   &        &   &   & d    & d    &&      &      \\
  &   &        &   &   & D & \ddots &   &  &&&&  &   &   &        &   &   & d    & d    &&      &      \\
  &   &        &   &   &   & \ddots & L &  &&&&  &   &   &        &   &   & \vdots    & \vdots    &&      &      \\
  &   &        &   &   &   &        & D &  &&&&  &   &   &        &   &   & d    & d    &&      &      \\
\hline
  &   &  &&&&  &   & L &   &        &   &        &   &   &        &   & D && d    & \ell &      &      \\
  &   &  &&&&  &   & D & L &        &   &        &   &   &        &   &   && d    & d    &      &      \\
  &   &  &&&&  &   &   & D & \ddots &   &        &   &   &        &   &   && d    & d    &      &      \\
  &   &  &&&&  &   &   &   & \ddots & L &        &   &   &        &   &   && \vdots    & \vdots    &      &      \\
  &   &  &&&&  &   &   &   &        & D &        &   &   &        &   &   && d    & d    &      &      \\
\hline
  &   &  &&&&  &   &   &   &        &   & \ddots &   &   &        & & \vdots && & \vdots &\ddots& \vdots \\
\hline
  &   &        &   &   &   &        &   &  &&&&  & L &   &        &   & D &&      &      &      & d    \\
  &   &        &   &   &   &        &   &  &&&&  & D & L &        &   &   &&      &      &      & d    \\
  &   &        &   &   &   &        &   &  &&&&  &   & D & \ddots &   &   &&      &      &      & d    \\
  &   &        &   &   &   &        &   &  &&&&  &   &   & \ddots & L &   &&      &      &      & \vdots    \\
  &   &        &   &   &   &        &   &  &&&&  &   &   &        & D &   &&      &      &      & d    
\end{array} 
\right)
   \label{mat:PretzelBAS}
\end{equation}

Observe that there are $k$ blocks, where the $i$-th block is $|n_i|\times(|n_i|-1)$, followed by some $1+(k-1)=k$ columns.  The first of these final columns has a single non-zero entry in each block:  a $D$ in the first position of the block, except for the first block which has an $L$ in the last position.  Of the remaining $k-1$ columns, the $i$-th has $|n_i|+|n_{i+1}|$ non-zero entries ($\ell d\ldots d$) appearing only in two consecutive blocks.

Either all the entries in a row have a bar or none of them do, and either all the rows in a block have barred entries or none of them do.  The first block in the matrix above has barred entries following $n_1<0$ in the pretzel knot depicted in the figure above.  In the \textit{enhanced activity matrix}, the writhe weighting $w(\varepsilon)$ will multiply each row by a monomial, and a Kasteleyn weighting given by Kauffman's trick $\kappa(\varepsilon)$ will involve signing some of the entries in the signed versions of these matrices.

Some care must be taken to order these activity letters correctly when assembling them into activity words.  The ordering of the edges in a perfect matching comes from the ordering of the original crossings, that is, by the ordering of the first vertex set or the rows of the matrix.

\begin{lemma}
\label{ActivityDimer}(Main Lemma)
Given a pretzel knot $P(n_1,n_2,\ldots,n_k)$, consider the usual diagram with the $n=|n_1|+\ldots+|n_k|$ crossings ordered from left to right and downward on the first column and then upward on the remaining columns.  Let $\Gamma$ be the balanced overlaid Tait graph for the diagram with the two omitted faces of the projection graph corresponding to the universal face and the upper deck supported by the columns as in Fig. \ref{fig:pretzel2TIKZ}.  

Summing over all perfect matchings $\mu$ in $\Gamma$ and taking the product over all edges $\varepsilon$ in the perfect matching,
\begin{equation}
\sum_{\mu}\prod_{\varepsilon\in\mu}\alpha(\varepsilon)=\sum_S a(S)
\end{equation}
gives the complete list of activity words $a(S)$ associated with spanning trees $S$ of the signed Tait graph $G$ associated with the diagram of the pretzel knot.
%
%
\end{lemma}

\begin{proof}
By Proposition \ref{prop:DimerTrees}, there is a bijection between the perfect matchings of the balanced overlaid Tait graph and the rooted spanning trees of the Tait graph associated with a knot diagram.  It is enough to show that the activity weighting of the perfect matching gives the activity word of the associated spanning tree.

Since there are $n=|n_1|+|n_2|+\ldots+|n_k|$ edges and $1+(|n_1|-1)+(|n_2|-1)+\ldots+(|n_k|-1)+1=n-k+2$ vertices in the Tait graph, each spanning tree $S$ omits exactly $k-1$ edges.  These $k-1$ edges must come from distinct columns in order for $S$ to be acyclic and connected.

The activity word associated with $S$ can be decomposed into the activity words of the $k$ columns, which are called 2-paths in this paper, and considered according to Lemma \ref{LdD}.

Suppose that an edge from each column except for the first is omitted.  Then by the lemma, the activity word of the first column will be ($L\ldots L$), and the activity words of the remaining columns will range from ($d D\ldots D$) to ($L\ldots L d D\ldots D$) to ($L\ldots L d$).

Suppose that an edge from each column except for the $i$-th is omitted for $i>1$.  Then by the lemma, the activity words of the first $i-1$ columns will range from ($\ell D\ldots D$) to ($L\ldots L d D\ldots D$) to ($L\ldots L d$), the activity word of the $i$-th column will be ($D \ldots D$), and the activity words of the remaining columns will range from ($d D\ldots D$) to ($L\ldots L d D\ldots D$) to ($L\ldots L d$).

As in the argument above, consider the location of the pivots in the last $k-1$ columns.  Note that if they both belong to the same $i$-th block, there must be a zero pivot in the $|n_i|\times(|n_i|-1)$ block to the left of it, and so this choice does not contribute to the permanent expansion.

When a pivot in each of the last $k-1$ columns is chosen, it forces the pivots of the corresponding blocks, with $L$'s chosen above the pivot row and $D$'s below in each block.  This then forces the pivot in the first of the last $k$ columns, which in turn forces the pivots in the remaining block to be uniform.

This gives the activity words described above.
\end{proof}

Evaluations of the activity words lead to the following main result.





\begin{theorem}
\label{DetJones}(Main Theorem)
Given a pretzel knot $K=P(n_1,\ldots,n_k)$, consider the usual diagram $D$ with the $n=|n_1|+\ldots+|n_k|$ crossings labelled from left to right and downward on the first column and then upward on the remaining columns.  Let $\Gamma$ be a balanced overlaid Tait graph for the diagram $D$ with the two omitted faces of the projection graph corresponding to the universal face and the upper deck supported by the columns as in Fig. \ref{fig:pretzel2TIKZ}.

Let $A$ be the signed enhanced activity matrix associated with $K$.  Then
\begin{equation}
\det(A)=V_K(t)
\end{equation}
gives the Jones polynomial $V_K(t)$ of $K$ up to sign.
\end{theorem}

\begin{proof}
Summing over all perfect matchings $\mu$ in $\Gamma$ and taking the product over all edges $\varepsilon$ in the perfect matching,
\begin{equation}
\sum_{\mu}\prod_{\varepsilon\in \mu}\alpha(\varepsilon)|_V=\langle K\rangle
\end{equation}
gives the Kauffman bracket polynomial $\langle K \rangle$ of $K$ up to sign, and
\begin{equation}
\sum_{\mu}\prod_{\varepsilon\in \mu}w(\varepsilon)\alpha(\varepsilon)|_V=V_K(t)
\end{equation}
gives the Jones polynomial $V_K(t)$ of $K$ up to sign.  This follows from Lemma \ref{ActivityDimer} by Theorem \ref{spJones} and the evaluations for the activity letters given in Table \ref{tab:ActivityEvaluations}.


Together with a Kasteleyn weighting $\kappa(\varepsilon)$, Proposition \ref{prop:DetDimer} gives the perfect matchings as the determinant of the matrix.
\end{proof}

\begin{remark}
One must use the ordering of the crossings given above to construct the matrix; there are orderings for which this result does not hold.
\end{remark}

For a further application using these activity words, see Section \ref{sec:Corollaries}.

\subsection{Beyond pretzel knots}
\label{subsec:beyond}

This is not the only class for which the activity weighting works.  Next is a technique to begin with a knot diagram that works and extend it outside of the class of pretzel knots.  Specifically it replaces the final crossing with any rational tangle (see for example \cite{KauffLamb:tangles} for more on this subject).

\begin{proposition}
\label{prop:subdividedouble}
(Subdivision/Doubling) 
Suppose the highest ordered edge $e_n$ in the Tait graph $G$ is incident with the omitted vertex and the omitted face.  Then if the activity weighting on a balanced overlaid Tait graph $\Gamma$ provides a dimer model for the knot diagram associated with the Tait graph $G$, this can be extended to the balanced overlaid Tait graph associated with a Tait graph $G\cup\{e_{n+1}\}$ that subdivides or doubles the edge $e_n$.
\end{proposition}

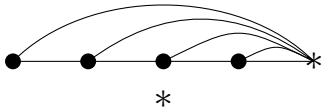
\begin{figure}[h]
\begin{center}

\begin{tikzpicture}
	\draw[-] (0,0) -- +(4,0);
	\draw[-] (0,0) .. controls (1,1) and (3,1).. (4,0);
	\draw[-] (1,0) .. controls (2,.75) and (3,.75) .. (4,0);
	\draw[-] (2,0) .. controls (3,.5) .. (4,0);
	\draw[-] (3,0) .. controls (3.5,.25) .. (4,0);
	\draw (4,0) node {\textbf{\Large $*$}};
	\fill[color=black] (0,0) circle (3pt);
	\fill[color=black] (1,0) circle (3pt);
	\fill[color=black] (2,0) circle (3pt);
	\fill[color=black] (3,0) circle (3pt);
	\draw (2,-.5) node {\textbf{\Large $*$}};
	
\end{tikzpicture}
	\caption{Doubling and subdivision leaves the class of pretzel knots.}
	\label{fig:SubdivideDouble}
\end{center}
\end{figure}

Note that for the ($n_1,\ldots,n_k$)-pretzel knot, if $|n_k|>1$ and the $n$-th edge is doubled, the end result is not a pretzel knot or link but a Montesinos knot or link.  See an example constructed from successive doubling and subdivisions on an edge in Fig. \ref{fig:SubdivideDouble}.

\begin{proof}
Consider the squared incidence matrix of the original knot diagram.  The row corresponding to the non-loop, non-bridge edge $e_n$, which is incident with both the omitted vertex and the omitted face, has only two non-zero entries.

If the edge is neither a bridge nor a loop, these entries are $D$ and $d$.  After subdividing this edge, the matrix gets a new row corresponding to the edge $e_{n+1}$ and a new column corresponding to the new vertex between $e_n$ and $e_{n+1}$.  The entries in this column are zero except for an $L$ and a $D$ in the $n$-th and ($n+1$)-st rows, respectively, and the last new entry is another $d$ in the ($n+1$)-st row below the first $d$ in the $n$-th row mentioned above.

Configurations in the determinant expansion for the final two terms in the new matrix have only three options: ($DD$) and ($dD$), which preserve all of the first $n$ choices of pivots, and ($Ld$), where the first $n-1$ choices are preserved and the $d$ of the $n$-th row gets replaced by the $d$ in the ($n+1$)-st row.  These are exactly the three possibilities for the activity words associated with the spanning trees by Lemma \ref{LdD}.

The dual case of doubling works similarly.
\end{proof}

\begin{remark}
\label{remark:Montesinos}
Note that a similar proof can be used to show that the first crossing can also be replaced with any rational tangle.  However, there are obstructions to this method preventing a crossing from the internal columns to be replaced.  Thus this construction cannot be used in general for the class of Montesinos knots (see for example treatment in Burde \cite{Bur}), whose ``columns'' are rational tangles instead of single twist regions.
\end{remark}

Although the Reidemeister moves do not alter the knot itself, they may significantly change the structure of the knot diagram.  Thus understanding how they affect the ordering of the crossings may lead to extending the activity matrix construction to larger knot classes.


The first Reidemeister move adds either bridges or loops to the Tait graph.  These are easily handled by the following:

\begin{proposition}
\label{prop:RMI}
(Reidemeister I) 
If the activity weighting on a balanced overlaid Tait graph $\Gamma$ associated with a knot diagram with signed Tait graph $G$ provides a dimer model, this can be extended to one whose Tait graph $G\cup \{e_{n+1}\}$ is the same as before together with an additional bridge or loop.
\end{proposition}

\begin{proof}
Since the edge $e_{n+1}$ is a bridge or a loop, this amounts only to adding the terms $L$ or $\ell$, respectively, to the end of the activity words, and this appears in the expansion because of a column with only a single non-zero entry.
\end{proof}

Each column of a pretzel knot contains copies of the same signed crossing, but the signed Tutte polynomial of more general signed graphs can also be considered.

\begin{proposition}
\label{property:RMII}
(Reidemeister II) 
If the activity weighting on a balanced overlaid Tait graph $\Gamma$ associated with a knot diagram with signed Tait graph $G$ provides a dimer model, this can be extended to one whose Tait graph $G\cup \{e_i',e_{i+1}'\}$ is the same as before together with two additional oppositely signed edges either in parallel or in series of an edge $e_i$ already in $G$.
\end{proposition}

\begin{proof}
This is a specialization of Lemma \ref{LdD} and its dual case Lemma \ref{ellDd}.
\end{proof}

Note that not all Reidemeister II moves are handled by this.  Propositions \ref{prop:RMI} and \ref{property:RMII} lead to the obvious question.

\begin{question}
Under what conditions can Reidemeister move III be used?
\end{question}

\section{Examples}
\label{sec:Examples}

\begin{example}  Fig. \ref{fig:trefoil} shows (A) the oriented trefoil depicted as the ($1,1,1$)-pretzel knot, (B) its corresponding (all-positive) Tait graph $G$, and (C) its corresponding balanced overlaid Tait graph $\Gamma$.  The thickened edges of the balanced overlaid Tait graph $\Gamma$ correspond to entries in the signed unenhanced activity matrix that are negatively signed by Kauffman's trick $\kappa(\varepsilon)$, giving the Kasteleyn weighting.  Note that the writhe of the diagram is $w(D)=-3$.

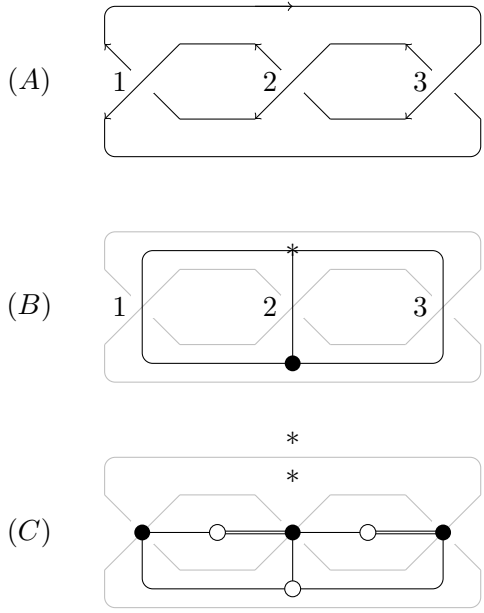
\begin{figure}
\begin{center}
\begin{tikzpicture}

\draw (-1,7) node {$(A)$};
\draw (-1,4) node {$(B)$};
\draw (-1,1) node {$(C)$};


\draw[-, rounded corners] (5,6.5) -- (5,6) -- (0,6) -- (0,6.5);

	\draw[<-] (0,7.5) -- +(1,-1);
	\fill[color=white] (.35,6.85) rectangle +(.3,.3);
	\draw[<-] (0,6.5) -- +(1,1);
\draw (.2,7) node {$1$};

\draw[-] (1,6.5) -- (2,6.5);
\draw[-] (1,7.5) -- (2,7.5);

	\draw[<-] (2,7.5) -- +(1,-1);
	\fill[color=white] (2.35,6.85) rectangle +(.3,.3);
	\draw[<-] (2,6.5) -- +(1,1);
\draw (2.2,7) node {$2$};

\draw[-] (3,6.5) -- (4,6.5);
\draw[-] (3,7.5) -- (4,7.5);

	\draw[<-] (4,7.5) -- +(1,-1);
	\fill[color=white] (4.35,6.85) rectangle +(.3,.3);
	\draw[<-] (4,6.5) -- +(1,1);
\draw (4.2,7) node {$3$};

\draw[-, rounded corners] (5,7.5) -- (5,8) -- (0,8) -- (0,7.5);
\draw[->] (2,8) -- (2.5,8);


\draw[-, rounded corners, lightgray] (5,3.5) -- (5,3) -- (0,3) -- (0,3.5);

	\draw[-, lightgray] (0,4.5) -- +(1,-1);
	\fill[color=white] (.35,3.85) rectangle +(.3,.3);
	\draw[-, lightgray] (0,3.5) -- +(1,1);
\draw (.2,4) node {$1$};

\draw[-, lightgray] (1,3.5) -- (2,3.5);
\draw[-, lightgray] (1,4.5) -- (2,4.5);

	\draw[-, lightgray] (2,4.5) -- +(1,-1);
	\fill[color=white] (2.35,3.85) rectangle +(.3,.3);
	\draw[-, lightgray] (2,3.5) -- +(1,1);
\draw (2.2,4) node {$2$};

\draw[-, lightgray] (3,3.5) -- (4,3.5);
\draw[-, lightgray] (3,4.5) -- (4,4.5);

	\draw[-, lightgray] (4,4.5) -- +(1,-1);
	\fill[color=white] (4.35,3.85) rectangle +(.3,.3);
	\draw[-, lightgray] (4,3.5) -- +(1,1);
\draw (4.2,4) node {$3$};

\draw[-, rounded corners, lightgray] (5,4.5) -- (5,5) -- (0,5) -- (0,4.5);

\draw[-, rounded corners] (.5,3.25) -- (.5,4.75) -- (4.5,4.75) -- (4.5,3.25) -- cycle;

	\draw (2.5,4.75) node {$*$};
	\fill[color=black] (2.5,3.25) circle (3pt);
\draw[-] (2.5,4.75) -- (2.5,3.25);



\draw[-, rounded corners, lightgray] (5,.5) -- (5,0) -- (0,0) -- (0,.5);

	\draw[-, lightgray] (0,1.5) -- +(1,-1);
	\fill[color=white] (.35,.85) rectangle +(.3,.3);
	\draw[-, lightgray] (0,.5) -- +(1,1);

\draw[-, lightgray] (1,.5) -- (2,.5);
\draw[-, lightgray] (1,1.5) -- (2,1.5);

	\draw[-, lightgray] (2,1.5) -- +(1,-1);
	\fill[color=white] (2.35,.85) rectangle +(.3,.3);
	\draw[-, lightgray] (2,.5) -- +(1,1);

\draw[-, lightgray] (3,.5) -- (4,.5);
\draw[-, lightgray] (3,1.5) -- (4,1.5);

	\draw[-, lightgray] (4,1.5) -- +(1,-1);
	\fill[color=white] (4.35,.85) rectangle +(.3,.3);
	\draw[-, lightgray] (4,.5) -- +(1,1);
	
\draw[-, rounded corners, lightgray] (5,1.5) -- (5,2) -- (0,2) -- (0,1.5);

\draw[-, rounded corners] (.5,1) -- (.5, .25) -- (4.5,.25) -- (4.5,1);


\draw[double] (1.5,1) -- (2.5,1);
\draw[double] (3.5,1) -- (4.5,1);


\draw[-] (2.5,1) -- (2.5,.25);
	\draw (2.5,2.25) node {$*$};
	\draw (2.5,1.75) node {$*$};
	\fill[color=white] (2.5,.25) circle (3pt);
	\draw (2.5,.25) circle (3pt);


	\fill[color=black] (.5,1) circle (3pt);
	\fill[color=black] (2.5,1) circle (3pt);
	\fill[color=black] (4.5,1) circle (3pt);


\draw[-] (.5,1) -- (1.5,1);
\draw[-] (2.5,1) -- (3.5,1);

	\fill[color=white] (1.5,1) circle (3pt);
	\draw (1.5,1) circle (3pt);
	\fill[color=white] (3.5,1) circle (3pt);
	\draw (3.5,1) circle (3pt);

\end{tikzpicture}
\end{center}
	\caption{(A) The (oriented) trefoil as the $(1,1,1)$-pretzel knot, (B) its corresponding Tait graph $G$, and (C) its corresponding balanced overlaid Tait graph $\Gamma$.}
	\label{fig:trefoil}
\end{figure}

The three spanning trees of the Tait graph are each of the three edges.  From Example \ref{example:trees}, the activity words are ($Ldd$), ($\ell Dd$), and ($\ell \ell D$), respectively.

The relatively trivial signed unenhanced activity matrix (without $w(\varepsilon)$, the writhe terms) has only the last columns and does not have any of the initial blocks:

\begin{equation} 
\left(
\begin{array}{c|cc}
L	& \ell	& 		 \\
\hline
D	& -d			& \ell	\\
\hline
D	& 			& -d		
\end{array} 
\right)
   \label{mat:example0}
\end{equation} 

One can see that the terms in the determinant expansion give the three activity words.  After evaluation, the determinant of this matrix is $-A^{-5}-A^3+A^7$.  Together with the term $(-A^{-3})^{-3}$ accounting for the writhe of the diagram, this gives $A^4+A^{12}-A^{16}=t^{-1}+t^{-3}-t^{-4}$, which is indeed the Jones polynomial of the trefoil.

\end{example}

\begin{example}
\label{example2}
Fig. \ref{fig:example1TIKZ} shows (A) the oriented ($-2,3,3$)-pretzel knot, also known as $8_{19}$, (B) its corresponding Tait graph $G$, and (C) its corresponding balanced overlaid Tait graph $\Gamma$.  The thickened edges of the signed Tait graph $G$ correspond to negative edges coming from negative crossings.  The thickened edges of the balanced overlaid Tait graph $\Gamma$ correspond to entries in the signed unenhanced activity matrix that are negatively signed by Kauffman's trick $\kappa(\varepsilon)$, giving the Kasteleyn weighting.  Note that the writhe of the diagram is $w(D)=8$.

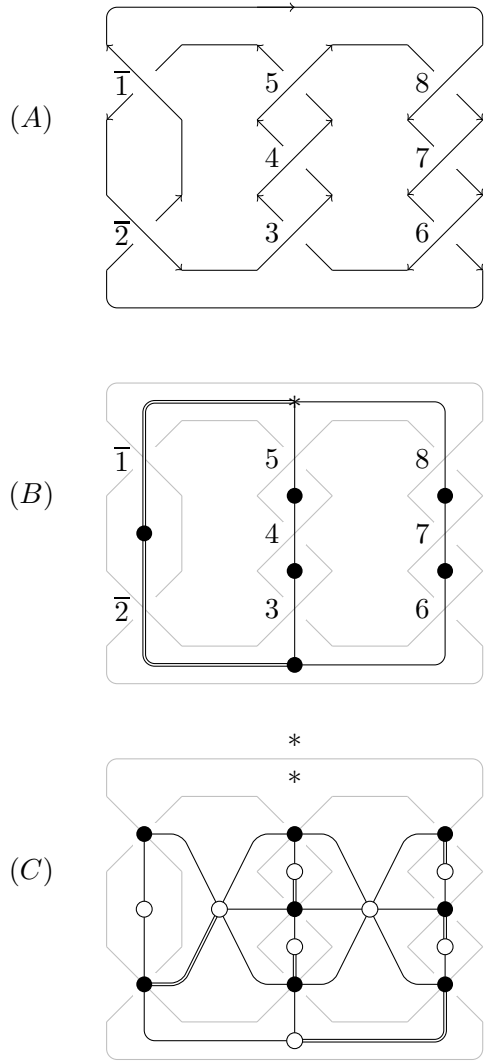
\begin{figure}
\begin{center}
\begin{tikzpicture}

\draw (-1,12.5) node {$(A)$};
\draw (-1,7.5) node {$(B)$};
\draw (-1,2.5) node {$(C)$};


\draw[-, rounded corners] [shift={(0,10)}] (5,.5) -- (5,0) -- (0,0) -- (0,.5);

	\draw[->] [shift={(0,10)}] (0,.5) -- +(1,1);
	\fill[color=white] [shift={(0,10)}] (.35,.85) rectangle +(.3,.3);
	\draw[->] [shift={(0,10)}] (0,1.5) -- +(1,-1);
\draw [shift={(0,10)}] (.2,1) node {$\overline{2}$};

\draw[-] [shift={(0,10)}] (0,1.5) -- (0,2.5);
\draw[-] [shift={(0,10)}] (1,1.5) -- (1,2.5);
	
	\draw[<-] [shift={(0,10)}] (0,2.5) -- +(1,1);
	\fill[color=white] [shift={(0,10)}] (.35,2.85) rectangle +(.3,.3);
	\draw[<-] [shift={(0,10)}] (0,3.5) -- +(1,-1);
\draw [shift={(0,10)}] (.2,3) node {$\overline{1}$};

\draw[-] [shift={(0,10)}] (1,.5) -- (2,.5);
\draw[-] [shift={(0,10)}] (1,3.5) -- (2,3.5);

	\draw[<-] [shift={(0,10)}] (2,1.5) -- +(1,-1);
	\fill[color=white] [shift={(0,10)}] (2.35,.85) rectangle +(.3,.3);
	\draw[->] [shift={(0,10)}] (2,.5) -- +(1,1);
\draw [shift={(0,10)}] (2.2,1) node {$3$};

	\draw[<-] [shift={(0,10)}] (2,2.5) -- +(1,-1);
	\fill[color=white] [shift={(0,10)}] (2.35,1.85) rectangle +(.3,.3);
	\draw[->] [shift={(0,10)}] (2,1.5) -- +(1,1);
\draw [shift={(0,10)}] (2.2,2) node {$4$};
	
	\draw[<-] [shift={(0,10)}] (2,3.5) -- +(1,-1);
	\fill[color=white] [shift={(0,10)}] (2.35,2.85) rectangle +(.3,.3);
	\draw[->] [shift={(0,10)}] (2,2.5) -- +(1,1);
\draw [shift={(0,10)}] (2.2,3) node {$5$};

\draw[-] [shift={(0,10)}] (3,.5) -- (4,.5);
\draw[-] [shift={(0,10)}] (3,3.5) -- (4,3.5);

	\draw[->] [shift={(0,10)}] (4,1.5) -- +(1,-1);
	\fill[color=white] [shift={(0,10)}] (4.35,.85) rectangle +(.3,.3);
	\draw[<-] [shift={(0,10)}] (4,.5) -- +(1,1);
\draw [shift={(0,10)}] (4.2,1) node {$6$};

	\draw[->] [shift={(0,10)}] (4,2.5) -- +(1,-1);
	\fill[color=white] [shift={(0,10)}] (4.35,1.85) rectangle +(.3,.3);
	\draw[<-] [shift={(0,10)}] (4,1.5) -- +(1,1);
\draw [shift={(0,10)}] (4.2,2) node {$7$};
	
	\draw[->] [shift={(0,10)}] (4,3.5) -- +(1,-1);
	\fill[color=white] [shift={(0,10)}] (4.35,2.85) rectangle +(.3,.3);
	\draw[<-] [shift={(0,10)}] (4,2.5) -- +(1,1);
\draw [shift={(0,10)}] (4.2,3) node {$8$};

\draw[-, rounded corners] [shift={(0,10)}] (5,3.5) -- (5,4) -- (0,4) -- (0,3.5);
\draw[->] [shift={(0,10)}] (2,4) -- (2.5,4);


\draw[-, rounded corners, lightgray] [shift={(0,5)}] (5,.5) -- (5,0) -- (0,0) -- (0,.5);

	\draw[-, lightgray] [shift={(0,5)}] (0,.5) -- +(1,1);
	\fill[color=white] [shift={(0,5)}] (.35,.85) rectangle +(.3,.3);
	\draw[-, lightgray] [shift={(0,5)}] (0,1.5) -- +(1,-1);
\draw [shift={(0,5)}] (.2,1) node {$\overline{2}$};

\draw[-, lightgray] [shift={(0,5)}] (0,1.5) -- (0,2.5);
\draw[-, lightgray] [shift={(0,5)}] (1,1.5) -- (1,2.5);
	
	\draw[-, lightgray] [shift={(0,5)}] (0,2.5) -- +(1,1);
	\fill[color=white] [shift={(0,5)}] (.35,2.85) rectangle +(.3,.3);
	\draw[-, lightgray] [shift={(0,5)}] (0,3.5) -- +(1,-1);
\draw [shift={(0,5)}] (.2,3) node {$\overline{1}$};

\draw[-, lightgray] [shift={(0,5)}] (1,.5) -- (2,.5);
\draw[-, lightgray] [shift={(0,5)}] (1,3.5) -- (2,3.5);

	\draw[-, lightgray] [shift={(0,5)}] (2,1.5) -- +(1,-1);
	\fill[color=white] [shift={(0,5)}] (2.35,.85) rectangle +(.3,.3);
	\draw[-, lightgray] [shift={(0,5)}] (2,.5) -- +(1,1);
\draw [shift={(0,5)}] (2.2,1) node {$3$};

	\draw[-, lightgray] [shift={(0,5)}] (2,2.5) -- +(1,-1);
	\fill[color=white] [shift={(0,5)}] (2.35,1.85) rectangle +(.3,.3);
	\draw[-, lightgray] [shift={(0,5)}] (2,1.5) -- +(1,1);
\draw [shift={(0,5)}] (2.2,2) node {$4$};
	
	\draw[-, lightgray] [shift={(0,5)}] (2,3.5) -- +(1,-1);
	\fill[color=white] [shift={(0,5)}] (2.35,2.85) rectangle +(.3,.3);
	\draw[-, lightgray] [shift={(0,5)}] (2,2.5) -- +(1,1);
\draw [shift={(0,5)}] (2.2,3) node {$5$};

\draw[-, lightgray] [shift={(0,5)}] (3,.5) -- (4,.5);
\draw[-, lightgray] [shift={(0,5)}] (3,3.5) -- (4,3.5);

	\draw[-, lightgray] [shift={(0,5)}] (4,1.5) -- +(1,-1);
	\fill[color=white] [shift={(0,5)}] (4.35,.85) rectangle +(.3,.3);
	\draw[-, lightgray] [shift={(0,5)}] (4,.5) -- +(1,1);
\draw [shift={(0,5)}] (4.2,1) node {$6$};

	\draw[-, lightgray] [shift={(0,5)}] (4,2.5) -- +(1,-1);
	\fill[color=white] [shift={(0,5)}] (4.35,1.85) rectangle +(.3,.3);
	\draw[-, lightgray] [shift={(0,5)}] (4,1.5) -- +(1,1);
\draw [shift={(0,5)}] (4.2,2) node {$7$};
	
	\draw[-, lightgray] [shift={(0,5)}] (4,3.5) -- +(1,-1);
	\fill[color=white] [shift={(0,5)}] (4.35,2.85) rectangle +(.3,.3);
	\draw[-, lightgray] [shift={(0,5)}] (4,2.5) -- +(1,1);
\draw [shift={(0,5)}] (4.2,3) node {$8$};

\draw[-, rounded corners, lightgray] [shift={(0,5)}] (5,3.5) -- (5,4) -- (0,4) -- (0,3.5);

\draw[-, rounded corners] [shift={(0,5)}] (2.5, .25) -- (4.5,.25) -- (4.5,3.75) -- (2.5,3.75);
\draw[-, rounded corners, double] [shift={(0,5)}] (2.5, .25) -- (.5,.25) -- (.5,3.75) -- (2.5,3.75);

	\fill[color=black] [shift={(0,5)}] (.5,2) circle (3pt);

	\draw [shift={(0,5)}] (2.5,3.75) node {$*$};
	\fill[color=black] [shift={(0,5)}] (2.5,2.5) circle (3pt);
	\fill[color=black] [shift={(0,5)}] (2.5,1.5) circle (3pt);
	\fill[color=black] [shift={(0,5)}] (2.5,.25) circle (3pt);
\draw[-] [shift={(0,5)}] (2.5,3.75) -- (2.5,.25);

	\fill[color=black] [shift={(0,5)}] (4.5,2.5) circle (3pt);
	\fill[color=black] [shift={(0,5)}] (4.5,1.5) circle (3pt);



\draw[-, rounded corners, lightgray] (5,.5) -- (5,0) -- (0,0) -- (0,.5);

	\draw[-, lightgray] (0,.5) -- +(1,1);
	\fill[color=white] (.35,.85) rectangle +(.3,.3);
	\draw[-, lightgray] (0,1.5) -- +(1,-1);

\draw[-, lightgray] (0,1.5) -- (0,2.5);
\draw[-, lightgray] (1,1.5) -- (1,2.5);
	
	\draw[-, lightgray] (0,2.5) -- +(1,1);
	\fill[color=white] (.35,2.85) rectangle +(.3,.3);
	\draw[-, lightgray] (0,3.5) -- +(1,-1);

\draw[-, lightgray] (1,.5) -- (2,.5);
\draw[-, lightgray] (1,3.5) -- (2,3.5);

	\draw[-, lightgray] (2,1.5) -- +(1,-1);
	\fill[color=white] (2.35,.85) rectangle +(.3,.3);
	\draw[-, lightgray] (2,.5) -- +(1,1);

	\draw[-, lightgray] (2,2.5) -- +(1,-1);
	\fill[color=white] (2.35,1.85) rectangle +(.3,.3);
	\draw[-, lightgray] (2,1.5) -- +(1,1);
	
	\draw[-, lightgray] (2,3.5) -- +(1,-1);
	\fill[color=white] (2.35,2.85) rectangle +(.3,.3);
	\draw[-, lightgray] (2,2.5) -- +(1,1);

\draw[-, lightgray] (3,.5) -- (4,.5);
\draw[-, lightgray] (3,3.5) -- (4,3.5);

	\draw[-, lightgray] (4,1.5) -- +(1,-1);
	\fill[color=white] (4.35,.85) rectangle +(.3,.3);
	\draw[-, lightgray] (4,.5) -- +(1,1);

	\draw[-, lightgray] (4,2.5) -- +(1,-1);
	\fill[color=white] (4.35,1.85) rectangle +(.3,.3);
	\draw[-, lightgray] (4,1.5) -- +(1,1);
	
	\draw[-, lightgray] (4,3.5) -- +(1,-1);
	\fill[color=white] (4.35,2.85) rectangle +(.3,.3);
	\draw[-, lightgray] (4,2.5) -- +(1,1);

\draw[-, rounded corners, lightgray] (5,3.5) -- (5,4) -- (0,4) -- (0,3.5);


\draw[rounded corners, double] (1.5,2) -- (1,1) -- (.5,1);

\draw[double] (2.5,1.5) -- (2.5,1);
\draw[double] (2.5,2.5) -- (2.5,2);

\draw[double] (4.5,3) -- (4.5,2.5);
\draw[double] (4.5,2) -- (4.5,1.5);
\draw[rounded corners, double] (4.5,1) -- (4.5,.25) -- (2.5,.25);


\draw[-, rounded corners] (.5,1) --(.5,.25) -- (2.5,.25);
\draw[-] (.5,3) -- (.5,1);
	\fill[color=white] (.5,2) circle (3pt);
	\draw (.5,2) circle (3pt);

\draw[-] (2.5,3) -- (2.5,2.5);
\draw[-] (2.5,2) -- (2.5,1.5);
\draw[-] (2.5,1) -- (2.5,.25);
	\draw (2.5,4.25) node {$*$};
	\draw (2.5,3.75) node {$*$};
	\fill[color=white] (2.5,2.5) circle (3pt);
	\draw (2.5,2.5) circle (3pt);
	\fill[color=white] (2.5,1.5) circle (3pt);
	\draw (2.5,1.5) circle (3pt);
	\fill[color=white] (2.5,.25) circle (3pt);
	\draw (2.5,.25) circle (3pt);

\draw[-] (4.5,2.5) -- (4.5,2);
\draw[-] (4.5,1.5) -- (4.5,1);
	\fill[color=white] (4.5,2.5) circle (3pt);
	\draw (4.5,2.5) circle (3pt);
	\fill[color=white] (4.5,1.5) circle (3pt);
	\draw (4.5,1.5) circle (3pt);


	\fill[color=black] (.5,3) circle (3pt);
	\fill[color=black] (.5,1) circle (3pt);

	\fill[color=black] (2.5,3) circle (3pt);
	\fill[color=black] (2.5,2) circle (3pt);
	\fill[color=black] (2.5,1) circle (3pt);

	\fill[color=black] (4.5,3) circle (3pt);
	\fill[color=black] (4.5,2) circle (3pt);
	\fill[color=black] (4.5,1) circle (3pt);


\draw[-] (1.5,2) -- (4.5,2);
\draw[-, rounded corners] (.5,3) -- (1,3) -- (1.5,2);
\draw[-, rounded corners] (2.5,3) -- (2,3) -- (1.5,2);
\draw[-, rounded corners] (1.5,2) -- (2,1) -- (2.5,1);
\draw[-, rounded corners] (2.5,3) -- (3,3) -- (3.5,2);
\draw[-, rounded corners] (3.5,2) -- (3,1) -- (2.5,1);
\draw[-, rounded corners] (4.5,3) -- (4,3) -- (3.5,2);
\draw[-, rounded corners] (3.5,2) -- (4,1) -- (4.5,1);

	\fill[color=white] (1.5,2) circle (3pt);
	\draw (1.5,2) circle (3pt);
	\fill[color=white] (3.5,2) circle (3pt);
	\draw (3.5,2) circle (3pt);

\end{tikzpicture}
\end{center}
	\caption{(A) The (oriented) ($-2,3,3$)-pretzel knot $8_{19}$, (B) its corresponding Tait graph $G$, and (C) its corresponding balanced overlaid Tait graph $\Gamma$.}
	\label{fig:example1TIKZ}
\end{figure}

The signed unenhanced activity matrix (without $w(\varepsilon)$, the writhe terms) is the following:

\begin{equation} 
\left(
\begin{array}{c|cc|cc|c|cc}
\overline{L} &   &   &   &   &              & \overline{\ell} &      \\
\overline{D} &   &   &   &   & \overline{L} & -\overline{d}    &      \\
\hline
             & -L &   &   &   &      D       &            d    & \ell \\
             & D & -L &   &   &              &            d    & d    \\
             &   & D &   &   &              &            d    & d    \\
\hline
             &   &   & L &   &      -D       &                 & d    \\
             &   &   & -D & L &              &                 & d    \\
             &   &   &   & -D &              &                 & d    
\end{array} 
\right)
   \label{mat:example1}
\end{equation} 

Taking the determinant after evaluation and together with the term $(-A^{-3})^{8}$ accounting for the writhe of the diagram, this gives $-A^{-32}+A^{-20}+A^{-12}=-t^8+t^5+t^3$, which is indeed the Jones polynomial of $8_{19}$.
\end{example}

\begin{example}  Given the ($-2,3,7$)-pretzel knot, which is useful in the construction of three manifolds, with crossings ordered downward on the first column, upward on the second, and then upward on the third, the signed unenhanced activity matrix (without $w(\varepsilon)$, the writhe terms) is the following:

\begin{equation} 
\left(
\begin{array}{c|cc|cccccc|c|cc}
\overline{L} &   &   &   &   &   &   &   &   &              & \overline{\ell} &      \\
\overline{D} &   &   &   &   &   &   &   &   & \overline{L} & -\overline{d}    &      \\
\hline
             & -L &   &   &   &   &   &   &   & D &        d    & \ell \\
             & D & -L &   &   &   &   &   &   &   &        d    & d    \\
             &   & D &   &   &   &   &   &   &   &        d    & d    \\
\hline
             &   &   & L &   &   &   &   &   & -D &             & d    \\
             &   &   & -D & L &   &   &   &   &   &             & d    \\
             &   &   &   & -D & L &   &   &   &   &             & d    \\
             &   &   &   &   & -D & L &   &   &   &             & d    \\
             &   &   &   &   &   & -D & L &   &   &             & d    \\
             &   &   &   &   &   &   & -D & L &   &             & d    \\
             &   &   &   &   &   &   &   & -D &   &             & d    \\
\end{array} 
\right)
   \label{mat:example2}
\end{equation} 

Taking the determinant after evaluation and together with the term $(-A^{-3})^{12}$ accounting for the writhe of the diagram, this gives $-A^{-40} + A^{-36} - A^{-32} + A^{-16} + A^{-8}=-t^{10}+t^9-t^8+t^2$.
\end{example}

It should be noted that although the negative signs in the matrices above can be easily obtained from the oriented knot diagram using Kauffman's trick, one can avoid this by computing the permanent of the unsigned matrix.

\section{Application to reduced Khovanov homology}
\label{sec:Corollaries}

Champanerkar and Kofman \cite{ChKo:mu} and \cite{ChKo:sp} define two gradings for a spanning tree $S$ of a Tait graph $G$ for a diagram $D$ with ordered crossings:
\begin{eqnarray*}
u(S)&=&\#L-\#\ell-\#\overline{L}+\#\overline{\ell}\\
v(S)&=&\#L+\#D
\end{eqnarray*}
and define a bigraded chain complex $\mathcal{C}(D)=\oplus_{u,v}\mathcal{C}^u_v(D)$, where
$$\mathcal{C}^u_v(D)=\mathbb{Z}\langle S\subset G|u(S)=u,v(S)=v\rangle,$$
to obtain $\mathcal{UC}(D)=\oplus_{u,v}(\mathcal{C}^u_v(D)+\overline{\mathcal{C}}^{u+2}_{v+1}(D))$, where $\overline{\mathcal{C}}^{u+2}_{v+1}(D)\cong\mathcal{C}^u_v(D)$ is a shifted copy.

Furthermore they prove that with a differential $\partial:C^u_v\rightarrow C^{u-1}_{v-1}$, this spanning tree complex  $\mathcal{UC}(D)$ is a deformation retract of the (unreduced) Khovanov complex.

Let $a(e,S)|_{Kh}$ denote the evaluation of the activity letters described in Table \ref{tab:KhActivityEvaluations} for spanning tree $S$ of the Tait graph $G$ of a knot diagram, and let $\alpha(\varepsilon)|_{Kh}$ be the corresponding evaluation in the balanced overlaid Tait graph $\Gamma$.

\begin{table}
	\centering
	\caption{Reduced Khovanov homology activity evaluations.}
{\begin{tabular}{|c|c|c|c|c|c|c|c|c|c|}
\hline
 & & & & & & & & & \\
$a(e,S)$ & activity letter of $e$ w.r.t. $S$ &	$L$ & $D$ & $\ell$ & $d$ & $\overline{L}$ & $\overline{D}$ & $\overline{\ell}$ & $\overline{d}$ \\
 & & & & & & & & & \\
\hline
 & & & & & & & & & \\
 & evaluation for the & & & & & & & & \\
$a(e,S)|_{Kh}$ & chain complex of reduced &	$uv$ & $v$ & $u^{-1}$ & $1$ & $u^{-1}$ & $1$ & $u$ & $1$ \\
 & Khovanov homology & & & & & & & & \\
 & & & & & & & & & \\
\hline
\end{tabular}}
	\label{tab:KhActivityEvaluations}
\end{table}

\begin{corollary}
\label{CKhDimer}
Given a pretzel knot $K=P(n_1,\ldots,n_k)$, consider the usual diagram $D$ with the $n=n_1+\ldots+n_k$ crossings labelled from left to right and downward on the first column and then upward on the remaining columns.  Let $\Gamma$ be a balanced overlaid Tait graph for the diagram $D$ with the two omitted faces of the projection graph corresponding to the universal face and the upper deck supported by the columns as in Fig. \ref{fig:pretzel2TIKZ}.

Summing over all perfect matchings $\mu$ in $\Gamma$ and taking the product over all edges $\varepsilon$ in the perfect matching,
\begin{equation}
\sum_{\mu}\prod_{\varepsilon\in \mu}\alpha(\varepsilon)|_{Kh}
\end{equation}
gives the two-variable Poincar\'e polynomial 
for the reduced Khovanov chain complex of $K$ up to sign.
\end{corollary}

\begin{proof}
This follows from Lemma \ref{ActivityDimer} by work of Champanerkar and Kofman and the evaluations for the activity letters given in Table \ref{tab:KhActivityEvaluations}.
\end{proof}

Champanerkar and Kofman distinguish between differentials coming from direct incidence and those induced from a collapse.  The former occurs if and only if the activity words of $S_1$ and $S_2$ differ by changing exactly two (not necessarily adjacent) letters in one of the following ways:
\begin{eqnarray*}
L\overline{d} & \rightarrow & d\overline{D}\\
\overline{d}D & \rightarrow & \overline{L}d\\
\overline{\ell}D & \rightarrow & \overline{D}d\\
D\overline{d} & \rightarrow & \ell\overline{D}.
\end{eqnarray*}

Remarkably these can easily be found in the (unsigned) unenhanced activity matrix.

\begin{corollary}
The differentials coming from direct incidence for the spanning tree model of reduced Khovanov homology corresponds to very particular $2\times2$ blocks of the (unsigned) unenhanced activity matrix (without the writhe and Kasteleyn weightings):  two rows and two columns who meet at four non-zero terms in each of the following configurations.
\begin{equation}
\left( \begin{array}{c|c}
L & d \\
\hline
\overline{D} & \overline{d} \end{array} \right)  \hspace{.2in}
\left( \begin{array}{c|c}
\overline{d} & \overline{L} \\
\hline
d & D \end{array} \right) \hspace{.2in}
\left( \begin{array}{c|c}
\overline{\ell} & \overline{D} \\
\hline
d & D \end{array} \right) \hspace{.2in}
\left( \begin{array}{c|c}
D & \ell \\
\hline
\overline{D} & \overline{d} \end{array} \right)
\end{equation}
\end{corollary}

Note that the order of the columns is not important; the second and third submatrices are written in this way to help illustrate the direction of the differential.

\begin{proof}
The subdeterminants of these $2\times2$ matrices give the direct incidence differentials because they only change two activity letters while keeping the rest of the activity word intact.
\end{proof}

Perhaps collections of edges that do not give spanning trees can be used to produce the differentials induced from a collapse through the context of the squared incidence matrix.  This would lead to a positive answer to the following question.

\begin{question}
Can the reduced Khovanov homology of pretzel knots can be computed via the (unsigned) unenhanced activity matrix alone?
\end{question}


\bibliographystyle{plain}
\bibliography{11FebruaryBibliography}

\end{document}